\documentclass[a4paper,10pt]{article}
\usepackage{appendix}
\usepackage{amsfonts}
\usepackage{epsfig}
\usepackage{amstext}
\usepackage{graphicx}
\usepackage{subfig}
\usepackage{amsthm}
\usepackage{listings}
\usepackage{verbatim}
\usepackage{amsmath,amssymb}
\usepackage{amsfonts}
\usepackage{graphicx}
\usepackage{pdfpages}
\usepackage{caption}
\usepackage[english]{babel}
\usepackage{hyperref}
\usepackage{float}
\usepackage[T1]{fontenc}
\usepackage[utf8]{inputenc}
\usepackage{authblk}
\usepackage{mathrsfs}

\usepackage{color}
\usepackage{hyperref}
\hypersetup{
    colorlinks,
    citecolor=black,
    filecolor=black,
    linkcolor=black,
    urlcolor=black
}
\numberwithin{equation}{section}
\numberwithin{figure}{section}

\begin{document}

\title{A boundary integral approach to linear and nonlinear transient wave scattering}
\author{Aihua Lin}
\author{Anastasiia Kuzmina}
\author{Per Kristen Jakobsen}
\affil{Department of Mathematics and Statistics, UIT the Arctic University of Norway, 9019 Troms\o, Norway}

\renewcommand\Authands{ and }
\maketitle

\begin{abstract}
In this paper we introduce a method for solving linear and nonlinear scattering problems for wave
equations using a new hybrid approach. This new approach consists of a  reformulation of the governing
 equations into a form that can be solved  by a  combination of a domain-based method and a boundary-integral method.
Our reformulation is aimed at a situation where we have a collection 
of compact scattering objects located in an otherwise homogeneous unbounded space.

The domain-based method is used to propagate the equations governing the wave field
 inside the scattering objects forward in time. The boundary integral
method is used to supply the domain-based method with the required boundary values
for the wave field.

  In this way the best features of both
methods come into play; the response inside the scattering objects, which
can be caused by both material  inhomogeneity and nonlinearities,
is easily taken into account using the domain-based method, and the boundary conditions supplied 
by the boundary integral method  makes it possible to confine the domain based method to the
inside of each scattering object.

\end{abstract}

\section{\bigskip Introduction}

Boundary integral formulations are well known in all areas of science and
technology and leads lead to highly efficient numerical algorithms for solving
partial differential equations. Their utility  are, in particular, evident for
scattering of waves from objects located in an unbounded space. For these
situations, one whole space dimension is taken out of the problem by reducing
the solution of the original PDEs to the solution of an integral equation
located on the boundaries of the scattering objects.

However, this reduction relies on the use of Green's functions and is therefore only
possible if the PDEs are linear. For computational reasons one is also  usually restricted to situations where
 the Green's functions are given by explicit formulas, and this rules out most situations where the materials are 
 inhomogeneous.  Since many problems of interest  involve scattering of waves  from objects that display both material  inhomogeneity and nonlinearity,
 boundary integral methods have appeared to be of
limited utility in computational science.  Adding, to the limited scope of the method, the fact that a somewhat
advanced mathematical machinery is needed to formulate PDEs in terms of
boundary integral equations, it is perhaps not hard to understand why the
method is not all that popular.

Domain-based methods, like the finite difference method and the finite element method,
 on the other hand, appear to be of
much wider utility. Their simple formulation and wide applicability to many types
 PDEs, both linear and nonlinear, have made them extremely popular in
the scientific computing community. In the context of scattering problems they do, however, have
problems of their own to contend with. These problems are of two quite distinct types. 

The first  type of problem is related to the fact that the scattering objects frequently represent abrupt changes
in material properties as compared to the properties of the surrounding homogeneous space. This abrupt
change leads to PDEs with discontinuous or near-discontinuous coefficients. Such features
are hard to represent accurately using finite element or finite difference methods. The favored
approach is to introduce multiple, interlinked grids, that are adjusted so that they conform to the
boundaries of the scattering objects. Generating such grids,  tailored to the possibly 
complex shape of the scattering objects, linking them together in a correct way and designing them in
such a way that the resulting numerical algorithm is accurate and stable, is challenging.  The approach has however  been refined over many years  and in general works quite well, but it certainly adds to the implementation complexity of these methods.

 The second type of problem is related to the fact that one can't grid the domain where the scattering objects are located
 for the simple reason that in almost all situations of interest this domain is unbounded. This problem is of course
 well known in the research community and the way it is resolved is to grid a computational box that is large enough
 to contain all scattering objects of interest. This can easily become a very large domain, leading to a very large number of degrees of freedom
 in the numerical algorithm. However, most of the time, the numerical algorithm associated with the domain has a simple structure
for which it is possible to design very fast implementations if the structure is taken into account. However, the introduction of the finite computational box in what is really a unbounded domain leads to the question of designing boundary conditions on the boundary of the box in such a way that it is fully transparent
to waves. This is not an easy thing to achieve, most approaches one can think of will in one way or another introduce an inhomogeneity
that will partly reflect waves hitting the boundary. This problem was first solved in a fully satisfactory manner for the case
of scattering of electromagnetic waves. The domain based method of choice for electromagnetic waves is the Finite Difference Time Domain method(FDTD) \cite{FDTD1},\cite{FDTD2},\cite{FDTD3}. This is, as the name indicates, a finite difference method, but a method  that has been designed to take into account the very special structure of Maxwell's equations.
The removal of reflections from the finite computational box was achieved by the introduction of a Perfectly Matched Layer(PML)\cite{PML1},\cite{PML2}. This 
amounts to adding a narrow layer of a specially constructed artificial material to the outside the computational box. The PML layer
is however only perfectly transparent to wave propagation if the grid has infinite resolution. For any finite grid there is still a small, but nonzero, reflection from
the boundary of the computational box. This can be reduced by making the PML layer thicker, but this leads to more degrees of freedom and thus an increasing computational
load. However, overall PML works well, and certainly much better than anything that came before it. There is no doubt that the introduction of PML was a breakthrough.

The use of PML was closely  linked to the special structure of Maxwell's equations. However,  it was soon realized that the same effect could be achieved by complexifying the
physical space outside the computational box and analytically continuing the fields into this complex spatial domain\cite{Chew},\cite{Johnson}.
 Significantly, this realization made the benefits of a  reflection less boundary condition available to all
kinds of scattering problems. However, the use of these reflection less boundary conditions certainly leads to an  increased computational load,  increased implementation complexity and also
to numerical stability issues that needs to be resolved.  It is at this point worth recalling that the boundary of the computational box is not part of the original physical
problem and that all the added implementation complexity and computational cost is spent trying to make it invisible after the choice of a domain method forced us to
put it there in the first place.

What we propose in this paper is to only apply the domain-based method inside
each scattering object. Firstly, this will reduce the size of the computational
grid enormously since we now only need to grid the inside of the scattering
objects. Secondly, our approach makes it possible to use different computational grids
for each scattering object, each grid tailored to the corresponding  object's  geometric shape,
without having to worry about the  inherent complexity caused by letting the different
grids meet up. Thirdly, it makes the introduction of a large computational
box, with its artificial boundary, redundant. In this way the computational
load is reduced by a large amount and we get rid of the implementation
complexity and instabilities associated with the boundary of the computational box.

However, the domain based method restricted to the inside of each scattering object  requires field values on the boundaries of
the scattering objects in order to be able to propagate the fields forward in
time. These boundary values will be supplied by a boundary integral method
derived from a space-time integral formulation of the PDEs one is seeking to solve. This
boundary integral method will take into account all the scattering and
 re-scattering of the solution to the PDEs in the unbounded domain outside the
scattering objects. Since the boundary integral method takes the radiation
condition at infinity into account explicitly, no finite computational box
with its artificial boundary conditions is needed.

This kind of idea for solving scattering problems was to our
knowledge first proposed in 1972 by Pattanayak and Wolf \cite{wolf} for the case of electromagnetic 
waves. They
discussed their ideas in the context of a generalization of the Ewald-Oseen
optical extinction theorem and we will because of this refer to our method as the
Ewald Oseen Scattering(EOS) formulation.

However, the paper by Pattanayak and Wolf only discussed stationary  linear scattering of
electromagnetic waves and they therefore did their integral formulation in
frequency space. This approach is not the right one when one is interested in transient scattering
from objects that are in general inhomogeneous and,  additionally, may have a  nonlinear response. What is needed for our approach is a space-time
integral formulation of the PDEs of interest.

In section 2 and 3 in this paper we illustrate our approach by  implementing our EOS formulation for two different 1D scattering problems. Both cases can be thought of as
toy models for the scattering of electromagnetic waves. This  should not be taken to mean that only models that in some way are related to electromagnetic scattering
can be subject to our approach, it merely reflects our particular interest in electromagnetic scattering. The way we see it, only one essential requirement needs to be fulfilled in order  for our method to be applicable;  it must be possible to 
derive an explicit integral formulation for the PDEs of interest. This means that at some point one needs to find the explicit expression for a Green's function for
some differential operator related to our PDEs. In general, it is difficult to find explicit expressions for Green's function belonging to nontrivial differential operators.  However the
 Green's function needed for our EOS formulation will always be of the infinite, homogeneous space type,  and explicit expressions for  such  Green's functions can frequently be found.
 
 The two models presented in this paper have been chosen for their simplicity, which makes them well suited for illustrating our EOS approach to scattering of waves. For more general and consequently more complicated cases there are really no new ideas beyond the technical details that must be mastered for each case  in order to derive the EOS formulation and implement it numerically. In order to explore the feasibility of our approach for more realistic and useful PDEs,  we have implemented our approach for several other cases, both 2d and 3D. In particular we have derived and implemented our EOS approach for the full 3D vector Maxwell's equations. The results of these investigations will be reported elsewhere.

 For both models we use an approach to testing stability and
accuracy of our implementations that involves what is known as \textit{artificial sources}. This method has probably been
around for a long time but apart from an application to the Navier-Stokes
equations \cite{Brehm11} we are not aware of any published work using this
method.  The method is based on the simple observation that if you add
arbitrary source terms to any system of PDEs then any function is a solution
for some choice of the source. Adding a source term typically introduce only 
trivial modifications to whichever numerical method used to solve the PDEs.
This essentially means that for any PDEs of interest, we can design particular
functions to test various critical aspects of the numerical method
related to numerical stability and accuracy.

This is a very simple approach to validating numerical implementations for PDEs
that deserves to be much better known than it apparently  is.

\bigskip

\section{The first scattering model; one way propagation}

Our first toy model, model 1, is%
\begin{align}
\varphi_{t}  &  =c_1\varphi_{x}+j,\label{A1}\nonumber\\
\rho_{t}  &  =-j_{x},\nonumber\\
j_{t}  &  =(\alpha-\beta\rho)\varphi-\gamma j\;\;\;\;a_0<x<a_1,
\end{align}
where  $\alpha,\beta$ and $\gamma$ are real parameters determining the ``material response'' part model 1 and where $\varphi=\varphi(x,t)$ is the ``electric field'', $j=j(x,t)$ the ``current
density'' and $\rho$ the ``charge density''. These quantities
are analogs for the corresponding quantities in Maxwell's equations. With this in mind, we observe that the second equation in the model (\ref{A1}) is a 1D version of the equation of
continuity from electromagnetics, and $c_1$ is the analog of the speed of light inside the ``material'' scattering object residing inside the interval $[a_0,a_1]$. The
charge density and current density are the material degrees of freedom and are therefore assumed to be confined to the
interval $[a_{0},a_{1}]$ on the real axis, whereas $\varphi$ is a field defined on the whole real axis. Thus the interval $[a_{0},a_{1}]$ is the
analog of a compact scattering object in the electromagnetic situation. Outside the scattering object the model
takes the form
\begin{align}
\varphi_{t}  &  =c_0\varphi_{x}+j_s\;\;\;\;x<a_0 \;\;\text{or}\;\;x>a_1,\label{A1.1}
\end{align}
where $c_0$ is the propagation speed for the electric field in the ``vacuum'' outside the scattering object and the
function $j_{s}(x,t)$ is a fixed source that has its support in a compact set
in the interval $x>a_{1}$. For the field $\varphi$ we impose the condition of continuity at the points $a_0$ and $a_1$.
The equation for the current density, $j$  is a radical simplification of a real current density model used to describe second harmonic generation in nonlinear optics \cite{Yong}.

\subsection{The EOS formulation}
In order to derive the EOS formulation for the
model (\ref{A1}), we will firstly need a space-time integral identity
involving the operator%
$$
L=\partial_{t}-v \partial_{x},
$$
where $v$ is some constant. Using integration by parts it is easy to see that
the following integral identity holds%
\begin{align}
&
{\displaystyle\int_{S\times T}}
dxdt\{L\varphi(x,t)\psi(x,t)-\varphi(x,t)L^{\dag}\psi(x,t)\}\label{A2}\nonumber\\
&  =%
{\displaystyle\int_{S}}
dx\varphi(x,t)\psi(x,t)|_{t_{0}}^{t_{1}}-v%
{\displaystyle\int_{T}}
dt\varphi(x,t)\psi(x,t)|_{x_{0}}^{x_{1}},
\end{align}
where $L^{\dag}=-\partial_{t}+v\partial_{x}$ is the formal adjoint of $L$ and
where $S=(x_{0},x_{1})$ and $T=(t_{0},t_{1})$ are open space and time intervals.

The second item we need in order to derive the EOS formulation for model
(\ref{A1}), is the advanced Green's function for the operator $L^{\dag}$. This
is a function $G=G(x,t,x^{\prime},t^{\prime})$ which is a solution to the
equation
\begin{equation}
L^{\dag}G(x,t,x^{\prime},t^{\prime})=\delta(t-t^{\prime})\delta(x-x^{\prime}),\nonumber
\end{equation}
and that vanishes when $t>t^{\prime}$. Using Fourier transforms we find that $G$
is given by%
\begin{equation}
G(x,t,x^{\prime},t^{\prime})=\theta(t^{\prime}-t)\delta(x^{\prime
}-x+v(t^{\prime}-t)), \label{A4}%
\end{equation}
where $\theta$ is the Heaviside step function with $\theta(x)=1$ for $x>0$
and zero otherwise.

We will now apply the integral identity (\ref{A2}) to each space interval
$(-\infty,a_{0})$, $(a_{0},a_{1})$ and $(a_{1},\infty)$. For the function $\psi$
we will substitute the advanced Green's function (\ref{A4}) and we will let $\varphi$ 
be the solution to equation (\ref{A1.1}) with vanishing initial condition, $\varphi(x,t_{0})=0$. We thus have a problem
where all solutions are purely source-generated.

  For the first interval, $(-\infty,a_{0})$, we let $\psi$ be the Green's
function%
\begin{equation}
G_0(x,t,x^{\prime},t^{\prime})\equiv\theta(t^{\prime}-t)\delta
(x^{\prime}-x+c_{0}(t^{\prime}-t)), \label{A6}%
\end{equation}
and $\varphi=\varphi_0$ be the solution to the
equation%
\begin{align}
\varphi_{0t}  &  =c_0\varphi_{0x},\label{A7}\nonumber\\
&  \Updownarrow\nonumber\\
L_0\varphi_0  &  =0.
\end{align}
Inserting (\ref{A6}), (\ref{A7}) and $S=(-\infty,a_{0})$ into the integral
identity (\ref{A2}), using the initial condition and the fact that the Green's
function is advanced, we get for $x$ in $(-\infty,a_{0})$%
\begin{align}
\varphi_0(x,t)  &  =c_{0}%
{\displaystyle\int_{t_{0}}^{t}}
dt^{\prime}\varphi(a_{0},t^{\prime})\delta(x-a_{0}+c_{0}(t-t^{\prime
}))\nonumber\\
&  -c_{0}\lim_{x^{\prime}\rightarrow-\infty}%
{\displaystyle\int_{t_{0}}^{t}}
dt^{\prime}\varphi_0(x^{\prime},t^{\prime})\delta(x-x^{\prime}%
+c_{0}(t-t^{\prime}))\nonumber\\
&  =c_{0}%
{\displaystyle\int_{t_{0}}^{t}}
dt^{\prime}\varphi_0(a_{0},t^{\prime})\delta(x-a_{0}+c_{0}(t-t^{\prime
})),\label{Identity1}
\end{align}
after interchanging the primed and unprimed variables. The last equality sign
follows because $x-x^{\prime}+c_{0}(t-t^{\prime})>0$ when $x^{\prime}<x$ for
all $t^{\prime}$ in the integration interval $(t_{0},t)$.

Note that when writing formula (\ref{Identity1}) we have made the substitution
\begin{equation*}
\varphi_{0}(a_0,\cdot)\equiv\lim_{x \rightarrow a_0^{-}} 
\varphi_{0}(x,\cdot). 
\end{equation*}
Similar substitutions will be made without comment later in this paper.

For the second interval, $(a_{0},a_{1})$, we let $\psi$ be the
Green's function%
\begin{equation}
G_1(x,t,x^{\prime},t^{\prime})\equiv\theta(t^{\prime}-t)\delta
(x^{\prime}-x+c_1(t^{\prime}-t)), \label{A9}%
\end{equation}
and $\varphi=\varphi_1$ be the solution to the
equation%
\begin{align}
\varphi_{1t}  &  =c_1\varphi_{1x}+j,\label{A10}\nonumber\\
&  \Updownarrow\nonumber\\
L_1\varphi  &  =j,
\end{align}
with vanishing initial conditions. Inserting (\ref{A9}),(\ref{A10}) and
$S=(a_{0},a_{1})$ into the integral identity (\ref{A2}), using the initial
condition and the fact that the Green's function is advanced, we get for $x$
in $(a_{0},a_{1})$
\begin{align}
\varphi_1(x,t)  &  =%
{\displaystyle\int_{a_{0}}^{a_{1}}}
dx^{\prime}%
{\displaystyle\int_{t_{0}}^{t}}
dt^{\prime}j(x^{\prime},t^{\prime})\delta(x-x^{\prime}+c_1(t-t^{\prime
}))\nonumber\\
&  +c_1%
{\displaystyle\int_{t_{0}}^{t}}
dt^{\prime}\varphi_1(a_{1},t^{\prime})\delta(x-a_{1}+c_1(t-t^{\prime
}))\nonumber\\
&  -c_1%
{\displaystyle\int_{t_{0}}^{t}}
dt^{\prime}\varphi_1(a_{0},t^{\prime})\delta(x-a_{0}+c_1(t-t^{\prime
}))\nonumber\\
&  =%
{\displaystyle\int_{a_{0}}^{a_{1}}}
dx^{\prime}%
{\displaystyle\int_{t_{0}}^{t}}
dt^{\prime}j(x^{\prime},t^{\prime})\delta(x-x^{\prime}+c_1(t-t^{\prime
}))\nonumber\\
&  +c_1%
{\displaystyle\int_{t_{0}}^{t}}
dt^{\prime}\varphi_1(a_{1},t^{\prime})\delta(x-a_{1}+c_1(t-t^{\prime
}))\label{Identity2},
\end{align}
after interchanging primed and unprimed variables. The last equality sign
follows because $x-a_{0}+c_1(t-t^{\prime})>0$ for all $t^{\prime}$ in the
integration interval when $a_{0}<x<a_{1\text{.}}$

 Finally, for the third integration interval, $(a_{1},\infty)$,  we let $\psi$ be
the Green's function%
\begin{equation}
G_0(x,t,x^{\prime},t^{\prime})\equiv\theta(t^{\prime}-t)\delta
(x^{\prime}-x+c_{0}(t^{\prime}-t)), \label{A12}%
\end{equation}
and $\varphi=\varphi_2$ be the solution to the
equation%
\begin{align}
\varphi_{2t}  &  =c_0\varphi_{2x}+j_{s},\label{A13}\nonumber\\
&  \Updownarrow\nonumber\\
L_0\varphi_2 &  =j_{s},
\end{align}
with vanishing initial conditions. Inserting (\ref{A12}),(\ref{A13}) and
$S=(a_{1},\infty)$ into the integral identity (\ref{A2}), using the initial
conditions and the fact that the Green's function is advanced, we get for $x$
in $(a_{1},\infty)$%
\begin{align}
\varphi_2(x,t)  &  =%
{\displaystyle\int_{a_{1}}^{\infty}}
dx^{\prime}%
{\displaystyle\int_{t_{0}}^{t}}
dt^{\prime}j_{s}(x^{\prime},t^{\prime})\delta(x-x^{\prime}+c_{0}(t-t^{\prime
}))\nonumber\\
&  +c_{0}\lim_{x^{\prime}\rightarrow\infty}%
{\displaystyle\int_{t_{0}}^{t}}
dt^{\prime}\varphi_2(x^{\prime},t^{\prime})\delta(x-x^{\prime}%
+c_{0}(t-t^{\prime}))\nonumber\\
&  -c_{0}%
{\displaystyle\int_{t_{0}}^{t}}
dt^{\prime}\varphi_2(a_{1},t^{\prime})\delta(x-a_{1}+c_{0}(t-t^{\prime
}))\nonumber\\
&  =%
{\displaystyle\int_{a_{1}}^{\infty}}
dx^{\prime}%
{\displaystyle\int_{t_{0}}^{t}}
dt^{\prime}j_{s}(x^{\prime},t^{\prime})\delta(x-x^{\prime}+c_{0}(t-t^{\prime
}))\label{Identity3},
\end{align}
after interchanging primed and unprimed variables. The third term vanishes
because $x-a_{1}+c_{0}(t-t^{\prime})>0$ for all $t^{\prime}$ in the
integration interval when $x>a_{1}$. The second term vanishes because
$x-x^{\prime}+c_{0}(t-t^{\prime})<0$ for all fixed $x>a_{1}$, $t>t_{0}\,$\ and
all $t^{\prime}$ in the integration interval $(t_{0},t)$ when $x^{\prime}$ is
large enough.

We now investigate the limit of these integral identities as $x$ approaches the
boundary points $\{a_{0},a_{1}\}$ of the open interval $(a_{0},a_{1})$ from
inside the interval and outside the interval. This will give us four equations
for the four quantities
\[
\varphi_0(a_{0},t),\varphi_1(a_{0},t),\varphi_1(a_{1},t),\varphi_2(a_{1},t).
\]
However, by assumption, acceptable solutions of model 1 are continuous
across the boundary points $\{a_{0},a_{1}\}.$  We therefore also have two
additional equations%
\begin{align*}
\varphi_0(a_{0},t)  &  =\varphi_1(a_{0},t),\nonumber\\
\varphi_1(a_{1},t)  &  =\varphi_2(a_{1},t).
\end{align*}
At this point we are faced with a problem; the four unknown quantities must satisfy 
 six linear equations. The problem is thus overdetermined and we
would not normally expect there to be any nontrivial solutions.

On the other hand, the equations, boundary conditions and source function $j_s$  that define model 1 do
determine a unique function $\varphi$. This function satisfies, by
construction, the integral identities (\ref{Identity1}),(\ref{Identity2}) and (\ref{Identity3}), whose limits yielded the overdetermined
system. So the overdetermined linear system does in fact have a solution.

There is a more direct way to see why the overdetermined system will have a solution.
Let us consider the inside of the scattering object, thus $x\in(a_0,a_1)$. Here, the field $\varphi$
is determined in terms of the current $j(x,t)$, and the boundary value $\varphi(a_1,t)$ by identity (\ref{Identity2})

\begin{align}
\varphi_1(x,t)  &  =%
{\displaystyle\int_{a_{0}}^{a_{1}}}
dx^{\prime}%
{\displaystyle\int_{t_{0}}^{t}}
dt^{\prime}j(x^{\prime},t^{\prime})\delta(x-x^{\prime}+c_1(t-t^{\prime
}))\nonumber\\
&  +c_1%
{\displaystyle\int_{t_{0}}^{t}}
dt^{\prime}\varphi_1(a_{1},t^{\prime})\delta(x-a_{1}+c_1(t-t^{\prime
})).\label{A11.1}
\end{align}
Naively, one would expect that we will get an equation determining the unknown boundary value
$\varphi(a_1,t)$, by taking the limit of (\ref{A11.1}) as $x$ approaches $a_1$ from below. However,
this would make the field inside the scattering object independent of the outside source, which from
a scattering point of view must be patently wrong. After all, it is the outside source $j_s(x,t)$ which
determines the field both outside and inside the scattering object. If this source is turned off the field
would simply be zero everywhere. So what is going on?

Note that if we actually take the limit of (\ref{A11.1}) we get the equation
\begin{equation}
0\;\varphi_1(a_1,t)=0,\nonumber
\end{equation}
which leaves the boundary value entirely arbitrary. If we analyze the rest of the overdetermined
 system in the same way, we find that one more equation for the boundary data is redundant, and
 that the two unknown boundary values, $\varphi_1(a_0,t)$ and $\varphi(a_1,t)$,  are uniquely determined by the following two equations
\begin{align}
\varphi_1(a_{0},t)  &  =%
{\displaystyle\int_{a_{0}}^{a_{1}}}
dx^{\prime}\theta(a_{0}-x^{\prime}+c_1(t-t_{0}))j(x^{\prime},t-\frac
{a_{1}-a_{0}}{c_1})\nonumber\\
&  +\theta(a_{0}-a_{1}+c_1(t-t_{0}))\varphi(a_{1},t-\frac{a_{1}-a_{0}}%
{c_1}),\label{A17}\\
\varphi_{1}(a_{1},t)  &  =\frac{1}{c_{0}}%
{\displaystyle\int_{a_{1}}^{\infty}}
dx^{\prime}\theta(a_{1}-x^{\prime}+c_{0}(t-t_{0}))j_{s}(x^{\prime
},t-\frac{x^{\prime}-a_{1}}{c_{0}}).\label{A18}
\end{align}
We emphasize the fact that we end up with an overdetermined system of linear
equations for the boundary values because this is a generic outcome when we derive
the EOS formulation for any given system of PDEs. We will see that  this very same problem will appear
when we discuss the second toy model in section three.

This problem has been recognized by the research community in the context of
space-time boundary integral formulation for the Maxwell's equations, and a
simple fix has been invented to resolve it.

However, as far as we know, the universal nature of this problem in the area
of space-time integral formulations of linear and nonlinear scattering
problems has not been recognized.

Observe that equation (\ref{A18}) determines the value of the field at the
boundary point $a_{1}$ in terms of the given external source $j_{s}$, and the equation
(\ref{A17}) determines the value of the field at the boundary point $a_{0}$
in terms of the current density  $j$ inside the scattering object and the field values at the boundary point $a_1$.

Equations (\ref{A1}) restricted to the the open interval $(a_{0},a_{1})$
together with the integral identities (\ref{A17}) and (\ref{A18}) define the
EOS formulation for model 1.

\subsection{Numerical implementation of the EOS formulation}\label{numericalf1}
In this section a numerical implementation of the EOS formulation for model 1 is presented. Many different numerical implementations are possible, the EOS formulation itself does not in any way dictate the use of some particular such implementation. However it does put some constraints on how we proceed with our method of choice.

    If our problem was to calculate the free-space propagation according to the first equation in (\ref{A1}) with vanishing $j$ the obvious choice would be to use the standard Lax-Wendroff method\cite{Lax} on a uniform space grid. However, the EOS formulation presents us with an integro-differential equation because the boundary update rule is defined in terms of integrals of the current density over the scattering domain $(a_0,a_1)$. Thus our grid must also give a good approximation to the integrals (\ref{A17}) and (\ref{A18}) which define the update rule. We will be looking for second order accuracy and would like to use the midpoint rule to approximate the integrals, and with this in mind, introduce the following nonuniform space grid  inside  the scattering object, $(a_0,a_1)$,
\begin{equation}
x_i=a_0+(i+0.5)\Delta x,\  i=0, 1, \cdots, N-1\label{OneFieldGrid},
\end{equation}
where $\Delta x =\frac{a_1-a_0}{N}$. The grid points  (\ref{OneFieldGrid}) will be called internal nodes in this paper. We also introduce a discrete time grid 
\begin{equation*}
t^n=n\Delta t,\  n=0, 1, \cdots. 
\end{equation*}
The values of the  parameter $\Delta t$ will, as usual, be bounded by the requirement of stability for the scheme. We will say a few words about this bound  later in the paper.

In order to get a numerical scheme of order two accuracy, we apply the Lax-Wendroff method to the first two equations of (\ref{A1}) and apply the modified Euler's method to the last equation of (\ref{A1}). Because of these choices the numerical scheme for iteration at the internal nodes takes the form
\begin{align}
\varphi_i^{n+1}
    =&\varphi_i^n+\Delta t \cdot ( c_1 \frac{\partial \varphi}{\partial x}+j )_i^n+\frac{1}{2} (\Delta t)^2( c_1^2 \frac{\partial^2 \varphi}{\partial x^2}+c_1 \frac{\partial j}{\partial x}+f)_i^n,\nonumber\\
\rho_i^{n+1}
    =&\rho_i^n+\Delta t \cdot ( - \frac{\partial j}{\partial x} )_i^n+\frac{1}{2} (\Delta t)^2( - \frac{\partial f}{\partial x})_i^n,\nonumber\\
 \bar j_i^{n+1}=&j_i^n+\Delta t \cdot f_i^n,\nonumber\\
  j_i^{n+1}=&\frac{1}{2}(j_i^n+\bar j_i^{n+1}+\Delta t\cdot  f(\rho_i^{n+1},\varphi_i^{n+1},\bar j_i^{n+1}))\label{LaxWendroff},
\end{align}
for $i=0, 1,\cdots, N$ and where $f=(\alpha-\beta\rho)\varphi-\gamma j$ . Except for the two internal nodes closest to the boundary points $a_0$ and $a_1$, the space derivatives are approximated to second order accuracy by the following standard finite difference formulas
\begin{align}
(\frac{\partial \phi}{\partial x})_i^n&=\frac{\phi_{i+1}^n-\phi_{i-1}^n}{2 \Delta x},\nonumber\\
(\frac{\partial^2 \phi}{\partial x^2})_i^n&=\frac{\phi_{i+1}^n-2 \phi_i^n+\phi_{i-1}^n}{(\Delta x)^2},\ \phi =\varphi,  j ,  f, \text{and}\  i=1,2,\cdots ,N-2.\label{spacem}
\end{align}
For the two internal nodes closest to the boundary, the standard, second order accurate, difference formulas, can not be used because the internal nodes are non-uniformly distributed in this part of the domain. 
For the field, $\varphi$, we must rather use the following second order accurate difference formulas for these two nodes
\begin{align}
(\frac{\partial \varphi}{\partial x})_0^n&=-\frac{1}{3\cdot \Delta x}(4\varphi_{a_0}^n-3\varphi_0^n-\varphi_1^n),\nonumber\\
(\frac{\partial^2 \varphi}{\partial x^2})_0^n&=\frac{4}{3\cdot (\Delta x)^2}(2\varphi_{a_0}^n-3\varphi_0^n+\varphi_1^n),\nonumber\\
(\frac{\partial \varphi}{\partial x})_{N-1}^n&=\frac{1}{3\cdot \Delta x}(4\varphi_{a_1}^n-3\varphi_{N-1}^n-\varphi_{N-2}^n),\nonumber\\
(\frac{\partial^2 \varphi}{\partial x^2})_{N-1}^n&=\frac{4}{3\cdot (\Delta x)^2}(2\varphi_{a_1}^n-3\varphi_{N-1}^n+\varphi_{N-2}^n).\label{spaceb}
\end{align}
The boundary value $\varphi_{a_0}^n$  needed in formulas (\ref{spaceb}) can be calculated from the discretized form of the integral update rules (\ref{A17})
\begin{equation}
\begin{split}
\varphi_{a_0}^{n+1}&=\frac{1}{c_1}\cdot \Delta x \cdot \sum_{i=0}^{N-1}\theta(t_{n+1}-t_0-\frac{x_i-a_0}{c_1})j(x_i,t_{n+1}-\frac{x_i-a_0}{c_1}),\\
&+\theta(t_{n+1}-t_0-\frac{a_1-a_0}{c_1})\varphi(a_1,t_{n+1}-\frac{a_1-a_0}{c_1})\label{BoundaryUpdate},
\end{split}
\end{equation}
while $\varphi_{a_1}^n$ is determined by the outside source using (\ref{A18}).

The current density, $j$, is entirely supported inside the scattering object and in general would be discontinuous at $a_0$ and $a_1$ if extended to the whole domain by making it zero external to the scattering object. Because of this, we need difference rules for $j$ at the nodes closest to the boundary points $a_0$ and $a_1$ that only depend on the values of $j$ on internal nodes. The following second order accurate difference rules for $j$ are of this type
\begin{align}
(\frac{\partial j}{\partial x})_0^n&=\frac{1}{2\Delta x}(4j_1^n-3j_0^n-j_2^n),\nonumber\\
(\frac{\partial j}{\partial x})_{N-1}^n&=-\frac{1}{2 \Delta x}(4j_{N-2}^n-3j_{N-1}^n-j_{N-3}^n).\label{jboundary}
\end{align}
It is evident that the discretized boundary update rule (\ref{BoundaryUpdate}) needs values of the current density that is located between the  grid points for the time grid $\{t^n\}$. This situation is general and will always arise when we seek numerical implementations of EOS formulations of PDEs. Some numerical interpolation scheme will always be needed to calculate the field values and/or the material variables between the time grid locations. Here we use a quadratic interpolation for values of the current density located between two time levels in order to maintain overall second order accuracy for our scheme. 

The iteration (\ref{LaxWendroff}) with the boundary update rule (\ref{BoundaryUpdate}) supplemented by the finite difference rules (\ref{spacem}),(\ref{spaceb}) and (\ref{jboundary}) constitute our numerical implementation of the EOS formulation for model 1.

\subsection{Artificial source test}
The basic idea behind the artificial source test, of some numerical scheme designed for a system of PDEs,  is to slightly modify the system by adding an arbitrary source to all the equations in the system. This modification typically lead to minimal modifications to the numerical scheme, where most of the effort and complexity are usually spent on the derivatives and nonlinear terms. For the equations, however, the presence of the sources  change the situation completely. This is because the presence of the added sources implies that {\it any} function is a solution to the equations for {\it some} choice of sources.

  With the risk of expanding on a perhaps already obvious idea; what we are saying is that if we have developed a numerical scheme for some system of differential equations $\mathscr{L}\psi=0$, we can with small modifications extend our scheme to the extended equation $\mathscr{L}\psi=g$ where $g$ is any given function. Given this, we test the numerical scheme by picking a function $\psi_0$, then use the equation to calculate the source function $g_0=\mathscr{L}\psi_0$ that ensure that our chosen function is a solution to the extended equation. Finally we run the numerical scheme with the calculated source function and find an approximate solution that we compare with the exact solution $\psi_0$.

 Mode 1 extended with artificial sources takes the form
\begin{align}
\varphi_{t}  &  =c_1 \varphi_{x}+j+g_1,\nonumber \\
\rho_{t}  &  =-j_{x}+g_2,\nonumber\\
j_{t}  &  =(\alpha-\beta\rho)\varphi-\gamma j+g_3, \label{A1s}
\end{align} where $g_1,$ $g_2,$ $g_3,$ are the artificial source functions. For some choice of functions $\varphi_0, j_0$ and  $\rho_0$ the corresponding source functions are computed by
\begin{equation*}
\begin{split}
g_{01}&=(\varphi_0)_t -c_1 (\varphi_0)_x -j_0, \\ \label{A1a}
g_{02}&=(\rho_0)_t +(j_0)_x ,\\
g_{03}&=(j_0)_t -(\alpha-\beta \rho_0)\varphi_0+\gamma j_0.
\end{split}
\end{equation*}
As our exact solution we choose 
\begin{align}
\varphi(x,t)=&\frac{2A_1}{\pi}\arctan(b^2t^2)e^{-\alpha_1(x-x_o+\beta_1(t-t_s))^2},\nonumber \\
j(x,t)=&A_2e^{-\frac{(x-x_j)^2}{\delta_1^2}-\frac{(t-t_j)^2}{\delta_2^2}},\nonumber\\
\rho(x,t)=&A_3e^{-\frac{(x-x_\rho)^2}{\delta_3^2}-\frac{(t-t_\rho)^2}{\delta_4^2}}\label{ExactSolution},
\end{align} 
which is nowhere near a solution to the equations  (\ref{A1.1}) defining  the unmodified model 1. Note that the chosen exact solution satisfies the vanishing of the initial data $\varphi(x,t_0=0)=0$, as it must in order to be consistent with the EOS formulation.

\begin{figure}[h!]
\centering
\includegraphics[width=12cm]{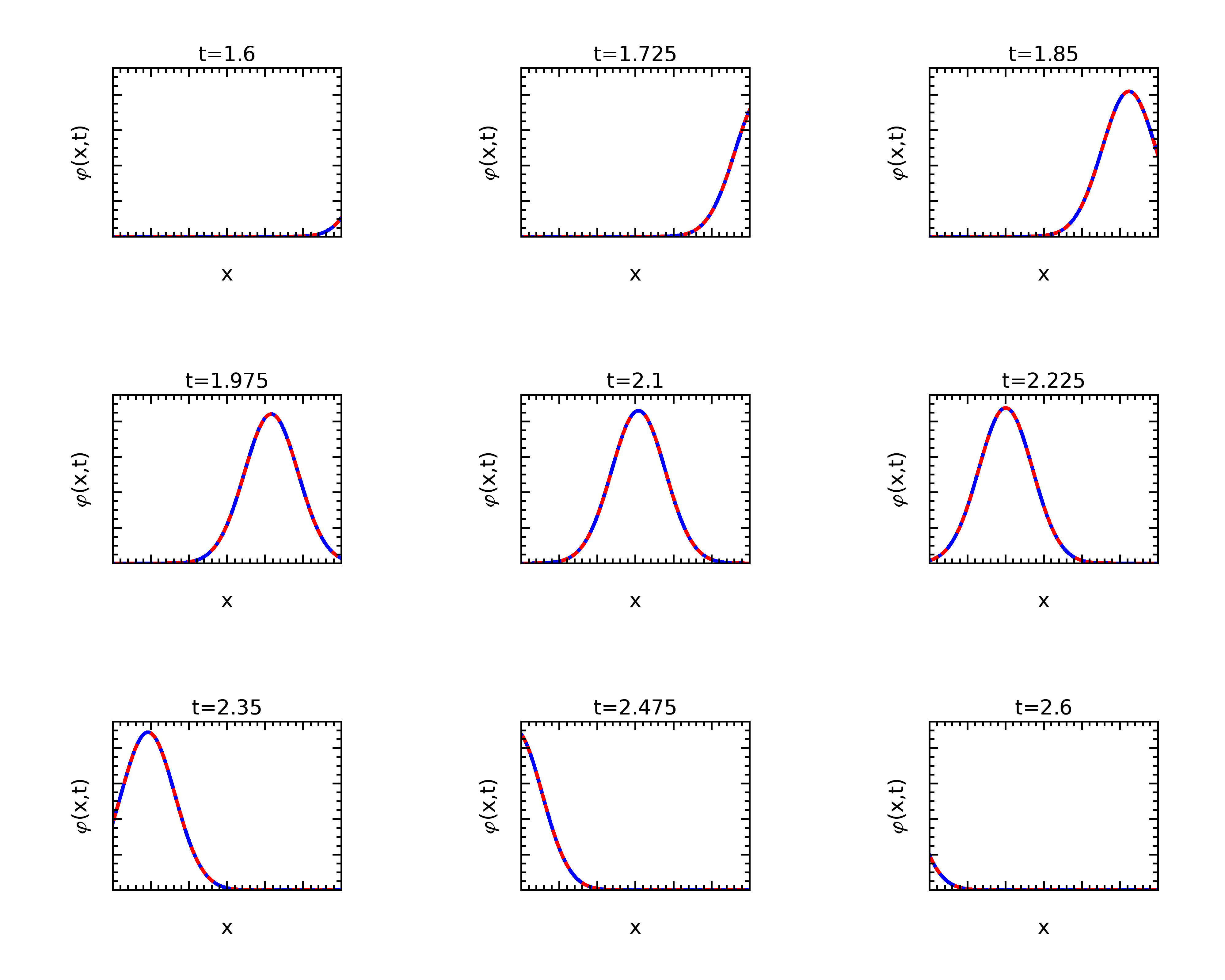}
\caption{Comparison between the numerical solution  and the exact solution for the source extended model 1 . Parameter values used are $a_0=0.0, a_1=3.0,
N=1600,$ $\alpha=-1.0,$ $\beta=0.3$, $\gamma=8.0$, $c=2.0$, $
c_0=1.0$, $A_1=1.0$,  $A_2=1.0$,  $A_3=1.0$, $b=1.0$, $\alpha_1=4.0, \beta_1=4.0, x_o=6.0, t_s=1.0$, $x_j=1.1, x_\rho=1.3, t_j=1.2, t_\rho=1.3, \delta_1=0.3, \delta_2=0.32, \delta_3=1.0, \delta_4=0.33 $. }
\label{fig:f1t}
\end{figure}

\begin{figure}[h!]
\centering
\includegraphics[width=12cm]{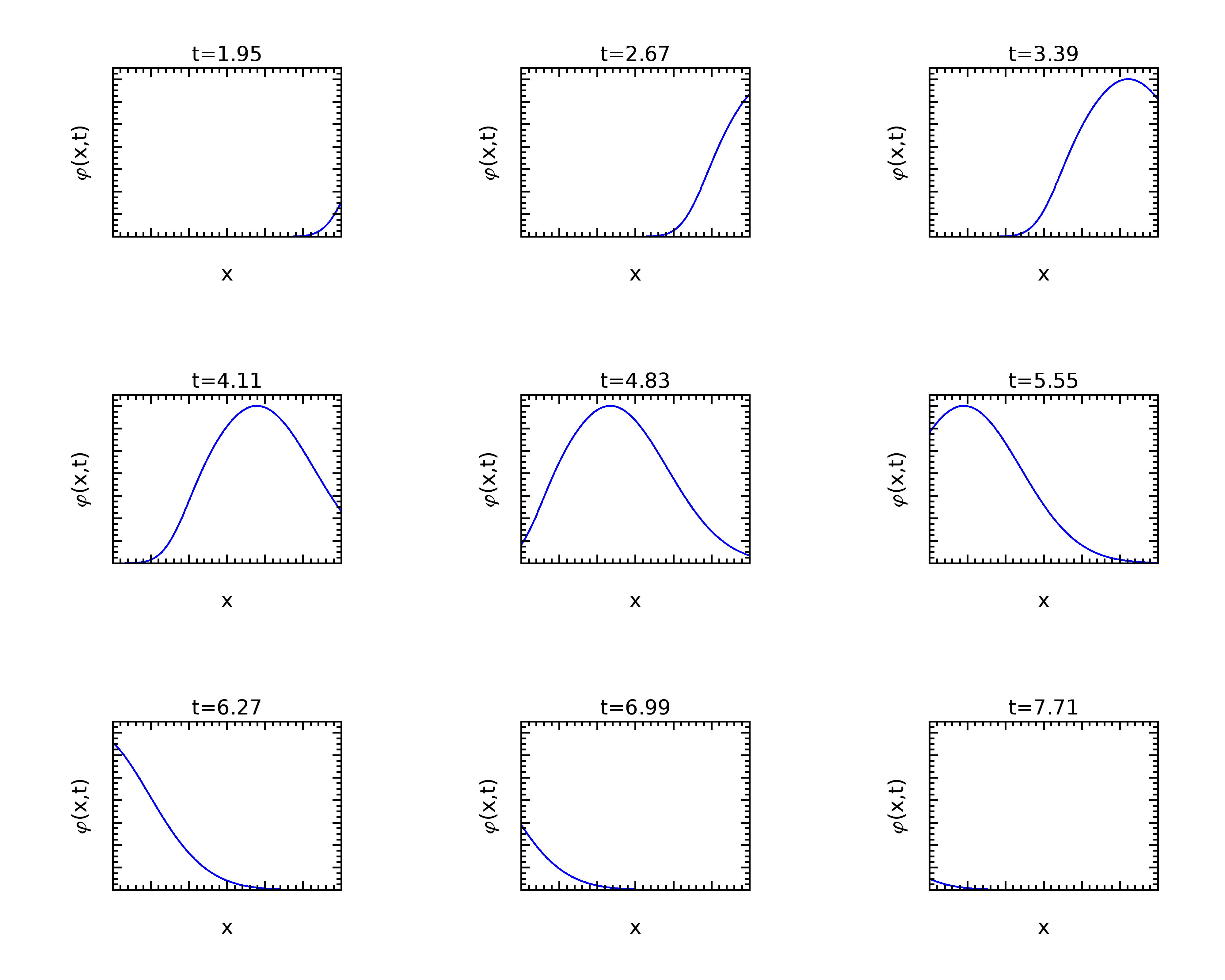}
\caption{A numerical solution of the EOS formulation for model 1 generated by an external source. The parameter values used are $a_0=0.0, a_1=3.0, N=1600, c=2.0,c_0=1.0, $ $\alpha=-1.0,$ $\beta=0.3$, $\gamma=8.0$. }
\label{fig:f1}
\end{figure}

The comparison between the exact solution (\ref{ExactSolution}) and the approximative solution generated by our numerical implementation of the EOS formulation of the source extended  model 1, (\ref{A1s}), is shown in figure \ref{fig:f1t} for some choice of the parameters. As we can see, the correspondence between the exact and approximative solution is excellent.
After having established that our implementation is accurate using the artificial source test, we show in figure  \ref{fig:f1} the numerical solution $\varphi$ of model 1, (\ref{A1}), where the system is driven by an outside source of the form
$$ j_s=5e^{-36(x-4)^2-4(t-0.5)^2},$$
 which is chosen so that no influence hit the boundary at $a_1$ before $t=0$. This will ensure that the initial condition $\varphi(x,t=0)=0$, underlying the EOS formulation of model 1, is satisfied.
 
   In these simulations we used a $\Delta t$ which is in the stable range for the numerical implementation, specifically we used $\Delta t =0.4 \frac{\Delta x}{c}$. Observe that the stability domain for our implementation of the EOS formulation is restricted compared to the stability domain for the underlying  Lax-Wendroff method on an infinite domain. The focus of the current paper is to derive the EOS formulation for two simple illustrative models and show that, using standard finite difference discretization of the EOS formulation, we get an accurate and stable representation of the solution to the scattering problems defined by the two toy models. A detailed and exhaustive discussion of the stability properties of the particular numerical scheme we have chosen for our implementation is less of a focus, and would be of interest only if such an analysis would lead to some kind of general statements with regard to numerical implementations of EOS formulations.  However, our experience with these numerical schemes, indicate that dimensionality is crucial with regards to stability and that the simplicity of the boundaries for 1D models makes them a poor guide to stability issues pertaining to EOS formulations in general.   Therefore, in order to stay focused on the main message of this paper, a discussion of the stability of our schemes for both toy models has been relegated to an appendix.

\section{The second scattering model; two way propagation}
Our second toy model, model 2  is%
\begin{align}
\varphi_{t}  &  =\mu_1\psi_{x}+j,\nonumber\\
\psi_{t}  &  =\nu_1\varphi_{x},\nonumber\\
\rho_{t}  &  =-j_{x},\nonumber\\
j_{t}  &  =(\alpha-\beta\rho)\varphi-\gamma j\label{B1.1}\;\;\;\;a_0<x<a_1,
\end{align}
where , like for model 1, $\varphi=\varphi(x,t)$, $j=j(x,t)$  and $\rho(x,t)$ are interpreted as ``electric field'', ``current density'' and ``charge density''. The additional field, $\psi(x,t)$ is interpreted as the ``magnetic'' field. The charge density and
current density will, as in model 1,  be confined to the interval $[a_{0},a_{1}]$
on the real axis whereas the fields $\varphi$ and $\psi$ are defined on the
whole real axis. The interval $[a_{0},a_{1}]$  is, like for model 1, the analog of
a compact scattering object in the electromagnetic situation. Outside the interval the model equations are 
\begin{align}
\varphi_{t}  &  =\mu_0\psi_{x}+j_{s},\nonumber\\
\psi_{t}  &  =\nu_0\varphi_{x}\label{B1.2},
\end{align}
where the function
$j_{s}(x,t)$ is a given source that, like for model 1, has its support on a compact set in the
interval $x>a_{1}$. The parameters $\mu_1,\mu_0,$
$\nu_1,\nu_0$ are "material" parameters. Using the translation $\mu\rightarrow\frac{1}{\epsilon}$ and $\nu\rightarrow\frac{1}{\mu}$ they are analogous for the electric permittivity, $\epsilon$, and the magnetic permeability, $\mu$,  inside and outside the scattering object.

\subsection{EOS formulation}
In order to derive the EOS formulation for model 2 (\ref{B1.1}),(\ref{B1.2}), we will firstly
need a space-time integral identity involving the matrix operator
\begin{equation}
L=\left(
\begin{array}
[c]{cc}%
\partial_{t} & -\mu\partial_{x}\\
-\nu\partial_{x} & \partial_{t}%
\end{array}
\right), \nonumber%
\end{equation}
where $\mu$ and $\nu$ are constants. The operator acts on vector valued
functions in the usual way
\begin{equation}
L\left(
\begin{array}
[c]{c}%
\varphi\\
\psi
\end{array}
\right)  =\left(
\begin{array}
[c]{c}%
\partial_t\varphi-\mu\partial_x\psi\\
\partial_t\psi-\nu\partial_x\varphi
\end{array}
\right). \nonumber%
\end{equation}
Using integration by parts, it is easy to derive the following integral
identity%
\begin{align}
&
{\displaystyle\int_{S\times T}}
dxdt\{AL\left(
\begin{array}
[c]{c}%
\varphi\\
\psi
\end{array}
\right)  (x,t)-L^{\dag}A\left(
\begin{array}
[c]{c}%
\varphi\\
\psi
\end{array}
\right)  (x,t)\}\label{B4}\nonumber\\
\nonumber\\
&  =%
{\displaystyle\int_{S}}
dxA\left(
\begin{array}
[c]{c}%
\varphi\\
\psi
\end{array}
\right)  (x,t)|_{t_{0}}^{t_{1}}+%
{\displaystyle\int_{T}}
dtB\left(
\begin{array}
[c]{c}%
\varphi\\
\psi
\end{array}
\right)  (x,t)|_{x_{0}}^{x_{1}},
\end{align}
where $S=(x_{0},x_{1})$ and $T=(t_{0},t_{1})$ are open space and time
intervals and where $\varphi$ and $\psi$ are smooth functions on the
space-time interval $S\times T$. Also $A=A(x,t)$ is a $2\times 2$ matrix valued function and $L^{\dag}$ is the formal adjoint to the operator $L$, and acts on the matrix
valued function $A$ in the following way%
\begin{equation}
L^{\dag}A=\left(
\begin{array}
[c]{cc}%
-\partial_{t}A_{11}+\nu\partial_{x}A_{12} & \mu\partial_{x}A_{11}-\partial
_{t}A_{12}\\
-\partial_{t}A_{21}+\nu\partial_{x}A_{22} & \mu\partial_{x}A_{21}-\partial
_{t}A_{22}%
\end{array}
\right).  \label{B6}%
\end{equation}
 $B$ is the $2\times2$ matrix valued function%
\begin{equation}
B=\left(
\begin{array}
[c]{cc}%
-\nu A_{12} & -\mu A_{11}\\
-\nu A_{22} & -\mu A_{21}%
\end{array}
\right).  \label{B5}%
\end{equation}
The second item we need in order to derive the EOS formulation for model
(\ref{B1.1}), (\ref{B1.2}), is the advanced Green's function for the operator $L^{\dag}$. This
is a $2\times2$ matrix valued function $G(x,t,x^{\prime},t^{\prime})$ that
satisfies the equation%
\begin{equation}
L^{\dag}G(x,t,x^{\prime},t^{\prime})=\delta(t-t^{\prime})\delta(x-x^{\prime})I,
\label{B8}%
\end{equation}
and that vanishes for $t>t^{\prime}$. In (\ref{B8}), $I$ is the $2\times2$
identity matrix

Using (\ref{B6}) in (\ref{B8}) we have the following system of four equations
for the components of $G$.%
\begin{align}
\partial_{t}G_{11}-\nu\partial_{x}G_{12}  &  =-\delta(t-t^{\prime}\nonumber%
)\delta(x-x^{\prime}),\label{B9}\\
\partial_{t}G_{12}-\mu\partial_{x}G_{11}  &  =0,\nonumber\\
\partial_{t}G_{21}-\nu\partial_{x}G_{22}  &  =0,\nonumber\\
\partial_{t}G_{22}-\mu\partial_{x}G_{21}  &  =-\delta(t-t^{\prime}%
)\delta(x-x^{\prime}).
\end{align}
We solve the system (\ref{B9}) using the Fourier transform and get%
\begin{align}
G(x,t,x^{\prime},t^{\prime})  &  =\frac{\theta(t^{\prime}-t)}{2c}\{\left(
\begin{array}
[c]{cc}%
c & \mu\\
\nu & c
\end{array}
\right)  \delta(x-x^{\prime}+c(t-t^{\prime}))\label{B10}\nonumber\\
\nonumber\\
&  +\left(
\begin{array}
[c]{cc}%
c & -\mu\\
-\nu & c
\end{array}
\right)  \delta(x-x^{\prime}-c(t-t^{\prime}))\},
\end{align}
where $c^{2}=\mu\nu$ and where $\theta(s)$ is the Heaviside step function.  Note that, using the identifications introduced while describing model 2 at the start of the current section,  the formula defining the speed, $c$,  is completely analogous to the one defining the speed of light in electromagnetics.

We will now apply the integral identity (\ref{B4}) to each space interval
$(-\infty,a_{0}),$ $(a_{0}, a_{1})$ and $(a_{1},\infty)$ with $A$ equal to the
advanced Green's function (\ref{B10}) for the corresponding interval and where $\varphi$ and $\psi$ are
solutions to the system (\ref{B1.1}),(\ref{B1.2}) with vanishing initial conditions $\varphi(x,t_{0})=\psi(x,t_{0})=0$.

For the first interval, $(-\infty,a_{0})$, we let $A$ be the Green's function%
\begin{align}
G_{0}(x,t,x^{\prime},t^{\prime})  &  =\frac{\theta(t^{\prime}-t)}{2c_{0}%
}\{\left(
\begin{array}
[c]{cc}%
c_{0} & \mu_0\\
\nu_0 & c_{0}%
\end{array}
\right)  \delta(x-x^{\prime}+c_{0}(t-t^{\prime}))\label{B12}\nonumber\\
\nonumber\\
&  +\left(
\begin{array}
[c]{cc}%
c_{0} & -\mu_0\\
-\nu_0 & c_{0}%
\end{array}
\right)  \delta(x-x^{\prime}-c_{0}(t-t^{\prime}))\},
\end{align}
where $c_0^{2}=\mu_0\nu_0$.
In this interval we let $\varphi=\varphi_{0},\psi=\psi_{0}$ be the solution to the system%
\begin{align}
\varphi_{0t}  &  =\mu_0\psi_{0x},\label{B13}\nonumber\\
\psi_{0t}  &  =\nu_0\varphi_{0x},\nonumber\\
&  \Updownarrow\nonumber\\
L_0\left(
\begin{array}
[c]{c}%
\varphi_{0}\\
\psi_{0}%
\end{array}
\right)   &  =0.
\end{align}
Inserting (\ref{B12}), (\ref{B13}) and $S=(-\infty,a_{0})$ into the integral
identity (\ref{B4}), using the initial conditions and the fact that the
Green's function is advanced, we get for $x$ in the interval $(-\infty,a_{0}%
)$. \
\begin{align}
\left(
\begin{array}
[c]{c}%
\varphi_{0}\\
\psi_{0}%
\end{array}
\right)  (x,t)  &  =-%
{\displaystyle\int_{t_{0}}^{t_{1}}}
dt^{\prime}B_{0}(a_{0},t^{\prime},x,t)\left(
\begin{array}
[c]{c}%
\varphi_{0}\\
\psi_{0}%
\end{array}
\right)  (a_{0},t^{\prime})\label{B14}\nonumber\\
&  +\lim_{R\rightarrow-\infty}%
{\displaystyle\int_{t_{0}}^{t_{1}}}
dt^{\prime}B_{0}(R,t^{\prime},x,t)\left(
\begin{array}
[c]{c}%
\varphi_{0}\\
\psi_{0}%
\end{array}
\right)  (R,t^{\prime}),
\end{align}
after interchanging primed and unprimed variables.

The function $B_{0}$ is from (\ref{B5})
\begin{align}
B_{0}(x^{\prime},t^{\prime},x,t)  &  =-\frac{\theta(t^{\prime}-t)}{2}\{\left(
\begin{array}
[c]{cc}%
c_{0} & \mu_0\\
\nu_0 & c_{0}%
\end{array}
\right)  \delta(x-x^{\prime}+c_{0}(t-t^{\prime}))\label{B15}\nonumber\\
\nonumber\\
&  +\left(
\begin{array}
[c]{cc}%
-c_{0} & \mu_0\\
\nu_0 & -c_{0}%
\end{array}
\right)  \delta(x-x^{\prime}-c_{0}(t-t^{\prime}))\}.
\end{align}
From (\ref{B15}) it is evident that the last term in (\ref{B14}) vanishes. This
is because for large enough $R$, the argument of the delta function does not
change sign in the interval of integration. Inserting the expression
(\ref{B15}) into (\ref{B14}) and changing to the variable defining the
argument of the delta function in the two integrals, we get that for $x$ in
$(-\infty,a_{0})$%
\begin{equation}
\left(
\begin{array}
[c]{c}%
\varphi_{0}\\
\psi_{0}%
\end{array}
\right)  (x,t)    =\frac{\theta(x-a_{0}+c_{0}(t-t_{0}))}{2c_{0}}\left(
\begin{array}
[c]{cc}%
c_{0} & \mu_0\\
\nu_0 & c_{0}%
\end{array}
\right)  \left(
\begin{array}
[c]{c}%
\varphi_{0}\\
\psi_{0}%
\end{array}
\right)  (a_{0},t+\frac{x-a_{0}}{c_{0}}).\label{B16}
\end{equation}
For the second interval, $(a_{0},a_{1})$, we let $A$ be the Green's function%
\begin{align}
G_{1}(x,t,x^{\prime},t^{\prime})  &  =\frac{\theta(t^{\prime}-t)}{2c_1%
}\{\left(
\begin{array}
[c]{cc}%
c _1& \mu_1\\
\nu_1 & c_1%
\end{array}
\right)  \delta(x-x^{\prime}+c_1(t-t^{\prime}))\label{B17}\nonumber\\
\nonumber\\
&  +\left(
\begin{array}
[c]{cc}%
c_1 & -\mu_1\\
-\nu_1 & c_1%
\end{array}
\right)  \delta(x-x^{\prime}-c_1(t-t^{\prime}))\}.
\end{align}
where $c_1^{2}=\mu_1\nu_1$.
In this interval,  the
functions $\varphi=\varphi_{1},\psi=\psi_{1}$ are the solutions to the system\ \ 
\begin{align}
\varphi_{1t}  &  =\mu_1\psi_{1x}+j,\label{B18}\nonumber \\
\psi_{1t}  &  =\nu_1\varphi_{1x},\nonumber\\
&  \Updownarrow\nonumber\\
L_1\left(
\begin{array}
[c]{c}%
\varphi_{1}\\
\psi_{1}%
\end{array}
\right)   &  =\left(
\begin{array}
[c]{c}%
j\\
0
\end{array}
\right). 
\end{align}
Inserting (\ref{B17}), (\ref{B18}) and $S=(a_{0},a_{1})$ in the integral
identity (\ref{B4}), using the vanishing  initial conditions and the fact that the
Green's function is advanced, we get for $x$ in the interval $(a_{0},a_{1})$.%
\begin{align}
\left(
\begin{array}
[c]{c}%
\varphi_{1}\\
\psi_{1}%
\end{array}
\right)  (x,t)  &  =%
{\displaystyle\int_{S\times T}}
dx^{\prime}dt^{\prime}G_{1}(x^{\prime},t^{\prime},x,t)\left(
\begin{array}
[c]{c}%
j\\
0
\end{array}
\right)  (x^{\prime},t^{\prime})\label{B19}\nonumber\\
&  -%
{\displaystyle\int_{t_{0}}^{t_{1}}}
dt^{\prime}B_{1}(a_{1},t^{\prime},x,t)\left(
\begin{array}
[c]{c}%
\varphi_{1}\\
\psi_{1}%
\end{array}
\right)  (a_{1},t)\nonumber\\
&  +%
{\displaystyle\int_{t_{0}}^{t_{1}}}
dt^{\prime}B_{1}(a_{0},t^{\prime},x,t)\left(
\begin{array}
[c]{c}%
\varphi_{1}\\
\psi_{1}%
\end{array}
\right)  (a_{0},t),
\end{align}
after interchanging primed and unprimed variables.

The function $B_{1}$ is from (\ref{B5})%
\begin{align}
B_{1}(x^{\prime},t^{\prime},x,t)  &  =-\frac{\theta(t^{\prime}-t)}{2}\{\left(
\begin{array}
[c]{cc}%
c_1 & \mu_1\\
\nu_1 & c_1%
\end{array}
\right)  \delta(x-x^{\prime}+c_1(t-t^{\prime}))\label{B20}\nonumber\\
\nonumber\\
&  +\left(
\begin{array}
[c]{cc}%
-c_1 & \mu_1\\
\nu_1 & -c_1%
\end{array}
\right)  \delta(x-x^{\prime}-c_1(t-t^{\prime}))\}.
\end{align}
Inserting (\ref{B17}) and (\ref{B20}) into (\ref{B19}), we get after changing
variables to the arguments in the delta functions that for $x$ in
$(a_{0},a_{1})$%
\begin{align}
&\left(
\begin{array}
[c]{c}%
\varphi_{1}\\
\psi_{1}%
\end{array}
\right)  (x,t)= \nonumber \\
\nonumber\\
 &\frac{1}{2c_1^{2}}\left(
\begin{array}
[c]{cc}%
c_1 & -\mu_1\\
-\nu_1 & c_1%
\end{array}
\right) \label{B21}
{\displaystyle\int_{a_{0}}^{x}}
dx^{\prime}\theta(c_1(t-t_{0})-(x-x^{\prime}))\binom{j}{0}(x^{\prime
},t-\frac{x-x^{\prime}}{c_1})\nonumber\\
\nonumber\\
&  +\frac{1}{2c_1^{2}}\left(
\begin{array}
[c]{cc}%
c_1 & \mu_1\\
\nu_1 & c_1%
\end{array}
\right)
{\displaystyle\int_{x}^{a_{1}}}
dx^{\prime}\theta(c_1(t-t_{0})-(x^{\prime}-x))\binom{j}{0}(x^{\prime
},t-\frac{x^{\prime}-x}{c_1})\nonumber\\
\nonumber\\
&  +\theta(c_1(t-t_{0})-(a_{1}-x))\frac{1}{2c_1}\left(
\begin{array}
[c]{cc}%
c_1 & \mu_1\\
\nu_1 & c_1%
\end{array}
\right) \left(
\begin{array}
[c]{c}%
\varphi_{1}\\
\psi_{1}%
\end{array}
\right)  (a_{1},t-\frac{a_{1}-x}{c_1})\nonumber\\
\nonumber\\
&  -\theta(c_1(t-t_{0})-(x-a_{0}))\frac{1}{2c_1}\left(
\begin{array}
[c]{cc}%
-c_1 & \mu_1\\
\nu_1 & -c_1%
\end{array}
\right)  \left(
\begin{array}
[c]{cc}%
\varphi_{1}\\
\psi_{1}%
\end{array}
\right)  (a_{0},t-\frac{x-a_{0}}{c_1}).
\end{align}
For the third interval, $(a_{1},\infty)$, we let $A$ be the Green's function%
\begin{align}
G_{0}(x,t,x^{\prime},t^{\prime})  &  =\frac{\theta(t^{\prime}-t)}{2c_{0}%
}\{\left(
\begin{array}
[c]{cc}%
c_{0} & \mu_0\\
\nu_0 & c_{0}%
\end{array}
\right)  \delta(x-x^{\prime}+c_{0}(t-t^{\prime}))\label{B22}\nonumber\\
\nonumber\\
&  +\left(
\begin{array}
[c]{cc}%
c_{0} & -\mu_0\\
-\nu_0 & c_{0}%
\end{array}
\right)  \delta(x-x^{\prime}-c_{0}(t-t^{\prime}))\}.
\end{align}
In this interval,  the
functions $\varphi=\varphi_{2},\psi=\psi_{2}$ are the solutions to the system%
\begin{align}
\varphi_{2t}  &  =\mu_0\psi_{2x}+j_{s},\label{B23}\nonumber\\
\psi_{2t}  &  =\nu_0\varphi_{2x},\nonumber\\
&  \Updownarrow\nonumber\\
L_0\left(
\begin{array}
[c]{c}%
\varphi_{2}\\
\psi_{2}%
\end{array}
\right)   &  =\left(
\begin{array}
[c]{c}%
j_{s}\\
0
\end{array}
\right). 
\end{align}
Inserting (\ref{B22}), (\ref{B23}) and $S=(a_{1},\infty)$ in the integral
identity (\ref{B4}), using the initial conditions and the fact that the
Green's function is advanced, we get for $x$ in the interval $(a_{1},\infty)$.%
\begin{align}
\left(
\begin{array}
[c]{c}%
\varphi_{2}\\
\psi_{2}%
\end{array}
\right)  (x,t)  &  =%
{\displaystyle\int_{S\times T}}
dx^{\prime}dt^{\prime}G_{0}(x^{\prime},t^{\prime},x,t)\left(
\begin{array}
[c]{c}%
j_{s}\\
0
\end{array}
\right)  (x^{\prime},t^{\prime})\label{B24}\nonumber\\
&  -\lim_{R\rightarrow\infty}%
{\displaystyle\int_{t_{0}}^{t_{1}}}
dt^{\prime}B_{0}(R,t^{\prime},x,t)\left(
\begin{array}
[c]{c}%
\varphi_{2}\\
\psi_{2}%
\end{array}
\right)  (R,t)\nonumber\\
&  +%
{\displaystyle\int_{t_{0}}^{t_{1}}}
dt^{\prime}B_{0}(a_{1},t^{\prime},x,t)\left(
\begin{array}
[c]{c}%
\varphi_{2}\\
\psi_{2}%
\end{array}
\right)  (a_{0},t),
\end{align}
after interchanging primed and unprimed variables.

Since the arguments of the delta functions in $B_{0}$ does not change sign in
the interval of integration, for $R$ big enough, it is clear that the second
term in (\ref{B24}) will vanish. Inserting (\ref{B22}) and (\ref{B15}) into
the remaining terms of (\ref{B24}), we get after changing variables to the
arguments in the delta functions that for $x$ in $(a_{1},\infty)$%
\begin{align}
&\left(
\begin{array}
[c]{c}%
\varphi_{2}\\
\psi_{2}%
\end{array}
\right)  (x,t)= \nonumber\\
\nonumber\\
&-\theta(c_{0}(t-t_{0})-(x-a_{1}))\frac{1}{2c_{0}}\left(
\begin{array}
[c]{cc}%
-c_{0} & \mu_0\\
\nu_0 & -c_{0}%
\end{array}
\right)\left(
\begin{array}
[c]{c}%
\varphi_{2}\\
\psi_{2}%
\end{array}
\right)  (a_{1},t-\frac{x-a_{1}}{c_{0}})\nonumber \\
&+\left(
\begin{array}
[c]{c}%
\varphi_{i}\\
\psi_{i}%
\end{array}
\right)  (x,t)\label{B26},
\end{align}
where $\varphi_{i}$ and $\psi_{i}$ are fields that are entirely determined by
the given source $j_{s}$%
\begin{align*}
&\left(
\begin{array}
[c]{c}%
\varphi_{i}\\
\psi_{i}%
\end{array}
\right)  (x,t)=  \\
\nonumber\\
&  \frac{1}{2c_{0}^{2}}\left(
\begin{array}
[c]{cc}%
c_{0} & -\mu_0\\
-\nu_0 & c_{0}%
\end{array}
\right)
{\displaystyle\int_{a_{1}}^{x}}
dx^{\prime}\theta(c_{0}(t-t_{0})-(x-x^{\prime}))\binom{j_{s}}{0}(x^{\prime
},t-\frac{x-x^{\prime}}{c_{0}})\nonumber\\
\nonumber\\
&  +\frac{1}{2c_{0}^{2}}\left(
\begin{array}
[c]{cc}%
c_{0} & \mu_0\\
\nu_0 & c_{0}%
\end{array}
\right) 
{\displaystyle\int_{x}^{\infty}}
dx^{\prime}\theta(c_{0}(t-t_{0})-(x^{\prime}-x))\binom{j_{s}}{0}(x^{\prime
},t-\frac{x^{\prime}-x}{c_{0}}).
\end{align*}
Taking the limit of the integral identities (\ref{B16}),(\ref{B21}) and
(\ref{B26}) as $x$ approaches the boundary points $\{a_{0},a_{1}\}$ from inside
and outside the interval $(a_{0},a_{1})$ we get
\begin{align}
&\left(
\begin{array}
[c]{cc}%
c_{0} & -\mu_0\\
-\nu_0 & c_{0}%
\end{array}
\right)  \left(
\begin{array}
[c]{c}%
\varphi_{0}\\
\psi_{0}%
\end{array}
\right)  (a_{0},t)=0, \label{B28}\\
\nonumber\\
&\left(
\begin{array}
[c]{cc}%
c_1 & \mu_1\\
\nu_1 & c_1%
\end{array}
\right)  \left(
\begin{array}
[c]{c}%
\varphi_{1}\\
\psi_{1}%
\end{array}
\right)  (a_{0},t)= \nonumber\\
\nonumber\\
&\frac{1}{c_1}\left(
\begin{array}
[c]{cc}%
c_1 & \mu_1\\
\nu_1 & c_1%
\end{array}
\right)
{\displaystyle\int_{a_{0}}^{a_{1}}}
dx^{\prime}\theta(c_1(t-t_{0})-(x^{\prime}-a_{0}))\binom{j}{0}(x^{\prime
},t-\frac{x^{\prime}-a_{0}}{c_1})\nonumber\\
\nonumber\\
&+\theta(c_1(t-t_{0})-(a_{1}-a_{0}))\left(
\begin{array}
[c]{cc}%
c_1 & \mu_1\\
\nu_1 & c_1%
\end{array}
\right)  \left(
\begin{array}
[c]{c}%
\varphi_{1}\\
\psi_{1}%
\end{array}
\right)  (a_{1},t-\frac{a_{1}-a_{0}}{c_1}),\label{B29}
\end{align}

\begin{align}
&\left(
\begin{array}
[c]{cc}%
c_1 & -\mu_1\\
-\nu_1 & c_1%
\end{array}
\right)  \left(
\begin{array}
[c]{c}%
\varphi_{1}\\
\psi_{1}%
\end{array}
\right)  (a_{1},t)=\nonumber\\
\nonumber\\
&
\frac{1}{c_1}\left(
\begin{array}
[c]{cc}%
c_1 & -\mu_1\\
-\nu_1 & c_1%
\end{array}
\right)
{\displaystyle\int_{a_{0}}^{a_{1}}}
dx^{\prime}\theta(c_1(t-t_{0})-(a_{1}-x^{\prime}))\binom{j}{0}(x^{\prime
},t-\frac{a_{1}-x^{\prime}}{c_1})\nonumber\\
\nonumber\\
&-\theta(c_1(t-t_{0})-(a_{1}-a_{0}))\left(
\begin{array}
[c]{cc}%
-c_1 & \mu_1\\
\nu_1 & -c_1%
\end{array}
\right)  \left(
\begin{array}
[c]{c}%
\varphi_{1}\\
\psi_{1}%
\end{array}
\right)  (a_{0},t-\frac{a_{1}-a_{0}}{c_1}),\label{B30}\\
\nonumber\\
&\left(
\begin{array}
[c]{cc}%
c_{0} & \mu_0\\
\nu_0 & c_{0}%
\end{array}
\right)  \left(
\begin{array}
[c]{c}%
\varphi_{2}\\
\psi_{2}%
\end{array}
\right)  (a_{1},t)=2c_{0}\left(
\begin{array}
[c]{c}%
\varphi_{i}\\
\psi_{i}%
\end{array}
\right)  (a_{1},t). \label{B31}%
\end{align}
Continuity of the fields at the boundary points $\{a_{0},a_{1}\}$, gives us
two additional equations,
\begin{align}
\left(
\begin{array}
[c]{c}%
\varphi_{0}\\
\psi_{0}%
\end{array}
\right)  (a_{0},t)  &  =\left(
\begin{array}
[c]{c}%
\varphi_{1}\\
\psi_{1}%
\end{array}
\right)  (a_{0},t),\label{B32}\\
\left(
\begin{array}
[c]{c}%
\varphi_{1}\\
\psi_{1}%
\end{array}
\right)  (a_{1},t)  &  =\left(
\begin{array}
[c]{c}%
\varphi_{2}\\
\psi_{2}%
\end{array}
\right)  (a_{1},t). \label{B33}%
\end{align}
Altogether we have six linear equations for the four vectors%
\[
\left(
\begin{array}
[c]{c}%
\varphi_{0}\\
\psi_{0}%
\end{array}
\right)  (a_{0},t),\left(
\begin{array}
[c]{c}%
\varphi_{1}\\
\psi_{1}%
\end{array}
\right)  (a_{0},t),\left(
\begin{array}
[c]{c}%
\varphi_{1}\\
\psi_{1}%
\end{array}
\right)  (a_{1},t),\left(
\begin{array}
[c]{c}%
\varphi_{2}\\
\psi_{2}%
\end{array}
\right)  (a_{1},t).
\]
Thus our system (\ref{B28})-(\ref{B33}) is overdetermined just like it was for
model 1. And just like for model 1, the system
(\ref{B28})-(\ref{B33}) contains equations that are redundant. Mathematically
 this is reflected in the fact that the determinant of the matrices%
\begin{equation}
\left(
\begin{array}
[c]{cc}%
c_{j} & \pm\mu_{j}\\
\pm\nu_{j} & c_{j}%
\end{array}
\right)  ,\text{ \ }j=0,1\;\;, \nonumber%
\end{equation}
are all zero. For the first toy model it was obvious which two equations were
redundant. Here it is not immediately clear which equations we can remove, and
this will also be the case if one write down the EOS formulation for more general systems of PDEs,
like for example Maxwell's equations.

For the system (\ref{B28})-(\ref{B33}), it is not very hard to identify the redundant
equations, but we will rather introduce a different approach that is in general
quite useful when working with the EOS formulations of PDEs.
This is the method that has been used by the research community that calculate
electromagnetic scattering from linear homogeneous scattering objects using a time dependent
integral formulation of Maxwell's equations. The {\it reason}
why this method has been used for the Maxwell's equations has not been clearly stated
 in the research literature. It has rather taken the form of a trick that is needed in order to
achieve stability and accuracy for the numerical implementation of the boundary
 formulation of electromagnetic scattering. 

The point is that even though the system (\ref{B28})-(\ref{B33}) is singular,
we know from it's construction that it has a solution which consists of the
boundary values coming from the unique solution to the system
(\ref{B1.1}),(\ref{B1.2}).

In terms of linear algebra, the situation is that for two given singular matrices $A$
and $B$, the system%
\begin{align}
A\mathbf{x}  &  =\mathbf{b}_{1},\label{B35}\nonumber\\
B\mathbf{x}  &  =\mathbf{b}_{2},
\end{align}
has a solution, $\mathbf{x}$. Let us assume that there are numbers $a$ and $b$
such that
\[
\det(aA+bB)\neq0.
\]
Given (\ref{B35}) it is clear that $\mathbf{x}$ is a solution to the linear
system%
\begin{equation}
(aA+bB)\mathbf{x}=\mathbf{b}_{1}+\mathbf{b}_{2}, \label{B36}%
\end{equation}
and since the system (\ref{B36}) is nonsingular, $\mathbf{x}$ is the unique
solution to the system. Finding numbers such that $aA+bB$ is nonsingular is in
general not difficult.

Let us apply this approach to the system (\ref{B28})-(\ref{B33}). Simply
adding together the equations give us a matrix%
\[
\left(
\begin{array}
[c]{cc}%
c_{0} & -\mu_0\\
-\nu_0 & c_{0}%
\end{array}
\right)  +\left(
\begin{array}
[c]{cc}%
c_1 & \mu_1\\
\nu_1 & c_1%
\end{array}
\right)  =\left(
\begin{array}
[c]{cc}%
c_1+c_{0} & \mu_1-\mu_0\\
\nu_1-\nu_0 & c_1+c_{0}%
\end{array}
\right), 
\] 
and
\[
\det\left(
\begin{array}
[c]{cc}%
c_1+c_{0} & \mu_1-\mu_0\\
\nu_1-\nu_0 & c_1+c_{0}%
\end{array}
\right)  =2c_1c_{0}+\mu_0\nu_1+\mu_1\nu_0, 
\]
which is nonzero since all the numbers $\nu_{i},\mu_{j},c_{j}$ are positive by
assumption. In a similar way, adding together (\ref{B30}) and (\ref{B31}) will
result in a nonsingular system. Thus from the singular system (\ref{B28}%
)-(\ref{B33}) we get the nonsingular system
\begin{align}
&\left(
\begin{array}
[c]{cc}%
c_1+c_{0} & \mu_1-\mu_0\\
\nu_1-\nu_0 & c_1+c_{0}%
\end{array}
\right)  \left(
\begin{array}
[c]{c}%
\varphi_{1}\\
\psi_{1}%
\end{array}
\right)  (a_{0},t)=\nonumber\\
\nonumber\\
&\frac{1}{c_1}\left(
\begin{array}
[c]{cc}%
c_1 & \mu_1\\
\nu_1 & c_1%
\end{array}
\right)
{\displaystyle\int_{a_{0}}^{a_{1}}}
dx^{\prime}\theta(c_1(t-t_{0})-(x^{\prime}-a_{0}))\binom{j}{0}(x^{\prime
},t-\frac{x^{\prime}-a_{0}}{c_1})\nonumber\\
\nonumber\\
&+\theta(c_1(t-t_{0})-(a_{1}-a_{0}))\left(
\begin{array}
[c]{cc}%
c_1 & \mu_1\\
\nu_1 & c_1%
\end{array}
\right)  \left(
\begin{array}
[c]{c}%
\varphi_{1}\\
\psi_{1}%
\end{array}
\right)  (a_{1},t-\frac{a_{1}-a_{0}}{c_1}),\label{B38.1}
\end{align}

\begin{align}
&\left(
\begin{array}
[c]{cc}%
c_{0}+c_1 & \mu_0-\mu_1\\
\nu_0-\nu_1 & c_1+c_{0}%
\end{array}
\right)  \left(
\begin{array}
[c]{c}%
\varphi_{1}\\
\psi_{1}%
\end{array}
\right)  (a_{1},t)=\nonumber\\
\nonumber\\
&\frac{1}{c_1}\left(
\begin{array}
[c]{cc}%
c_1 & -\mu_1\\
-\nu_1 & c_1%
\end{array}
\right)
{\displaystyle\int_{a_{0}}^{a_{1}}}
dx^{\prime}\theta(c_1(t-t_{0})-(a_{1}-x^{\prime}))\binom{j}{0}(x^{\prime
},t-\frac{a_{1}-x^{\prime}}{c_1})\nonumber\\
\nonumber\\
&-\theta(c_1(t-t_{0})-(a_{1}-a_{0}))\left(
\begin{array}
[c]{cc}%
-c_1 & \mu_1\\
\nu_1 & -c_1%
\end{array}
\right)  \left(
\begin{array}
[c]{c}%
\varphi_{1}\\
\psi_{1}%
\end{array}
\right)  (a_{0},t-\frac{a_{1}-a_{0}}{c_1})\nonumber\\
&+2c_{0}\left(
\begin{array}
[c]{c}%
\varphi_{i}\\
\psi_{i}%
\end{array}
\right)  (a_{1},t)\label{B38.2}.
\end{align}
The system (\ref{B38.1}),(\ref{B38.2}), which determine the boundary values of the
fields in term of  internal and external current densities, together with the differential
equations (\ref{B1.1}), restricted to the  inside the scattering object $(a_0,a_1)$, constitute the EOS formulation for model 2.

Our numerical implementation of model 2 is very similar to the one for model 1 and we therefore relegate the details of the implementation and artificial source test to an appendix. The results we get from this implementation is summed up in the next two figures.

\begin{figure}[h!]
\centering
\includegraphics[height=12cm, width=12cm]{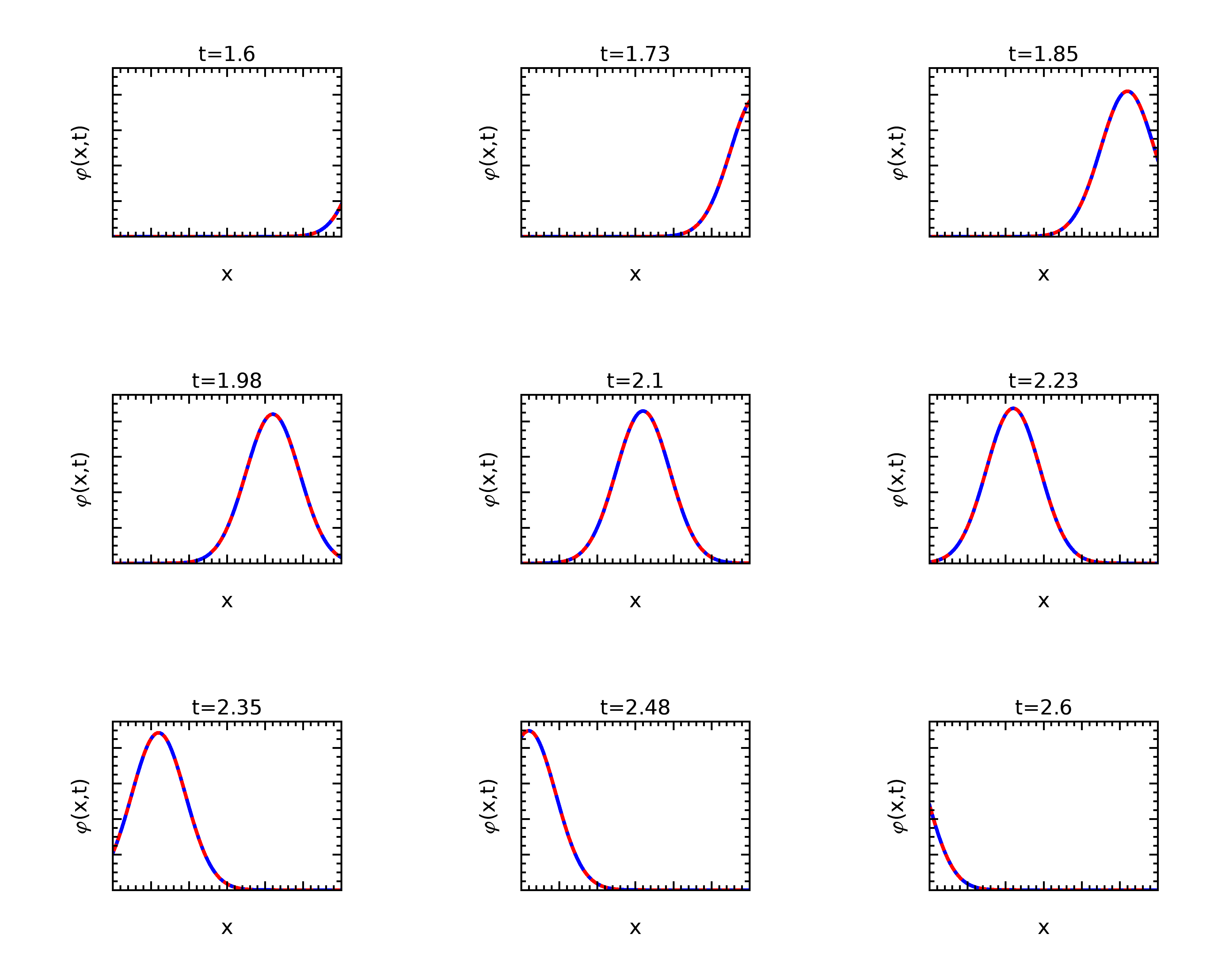}
\caption{Comparison between the numerical solution  and the exact solution for the source extended model 2 . Parameter values used are  $a_0=0.0, a_1=3.0, N=1600,$ $\alpha=-1.0,$ $\beta=0.3$, $\gamma=8.0$, $\mu=2.0, \nu=2.0$, $
\mu_0=1.0, \nu_0=1.0$, $A_1=1.0$,  $A_2=1.0$,  $A_3=1.0$, $A_4=1.0$, $b_1=1.0$, $b_2=1.0$, $\alpha_1=4.0, \beta_1=4.0, \alpha_2=4.0, \beta_2=4.0, x_o=6.0, t_s=1.0$, $x_j=1.1, x_\rho=1.3, t_j=1.2, t_\rho=1.3, \delta_1=0.3, \delta_2=0.32, \delta_3=1.0, \delta_4=0.33 $.}
\label{SourceTestModel2}
\end{figure}

In figure \ref{SourceTestModel2} we compare the numerical and exact solution of the EOS formulation for the source extended  model 2. The exact solution we used for this test  is 
\begin{equation*}
\begin{split}
\varphi(x,t)&=\frac{2A_1}{\pi}\arctan(b_1^2t^2)e^{-\alpha_1(x-x_o+\beta_1(t-t_s))^2},\\
\psi(x,t)&=\frac{2A_2}{\pi}\arctan(b_2^2t^2)e^{-\alpha_2(x-x_o+\beta_2(t-t_s))^2},\\
j(x,t)&=A_3e^{-\frac{(x-x_j)^2}{\delta_1^2}-\frac{(t-t_j)^2}{\delta_2^2}},\\
\rho(x,t)&=A_4e^{-\frac{(x-x_\rho)^2}{\delta_3^2}-\frac{(t-t_\rho)^2}{\delta_4^2}}.
\end{split}
\end{equation*}
Our implementation clearly passes the artificial source test with flying colors.

\begin{figure}[h!]
\centering
\includegraphics[height=12cm, width=12cm]{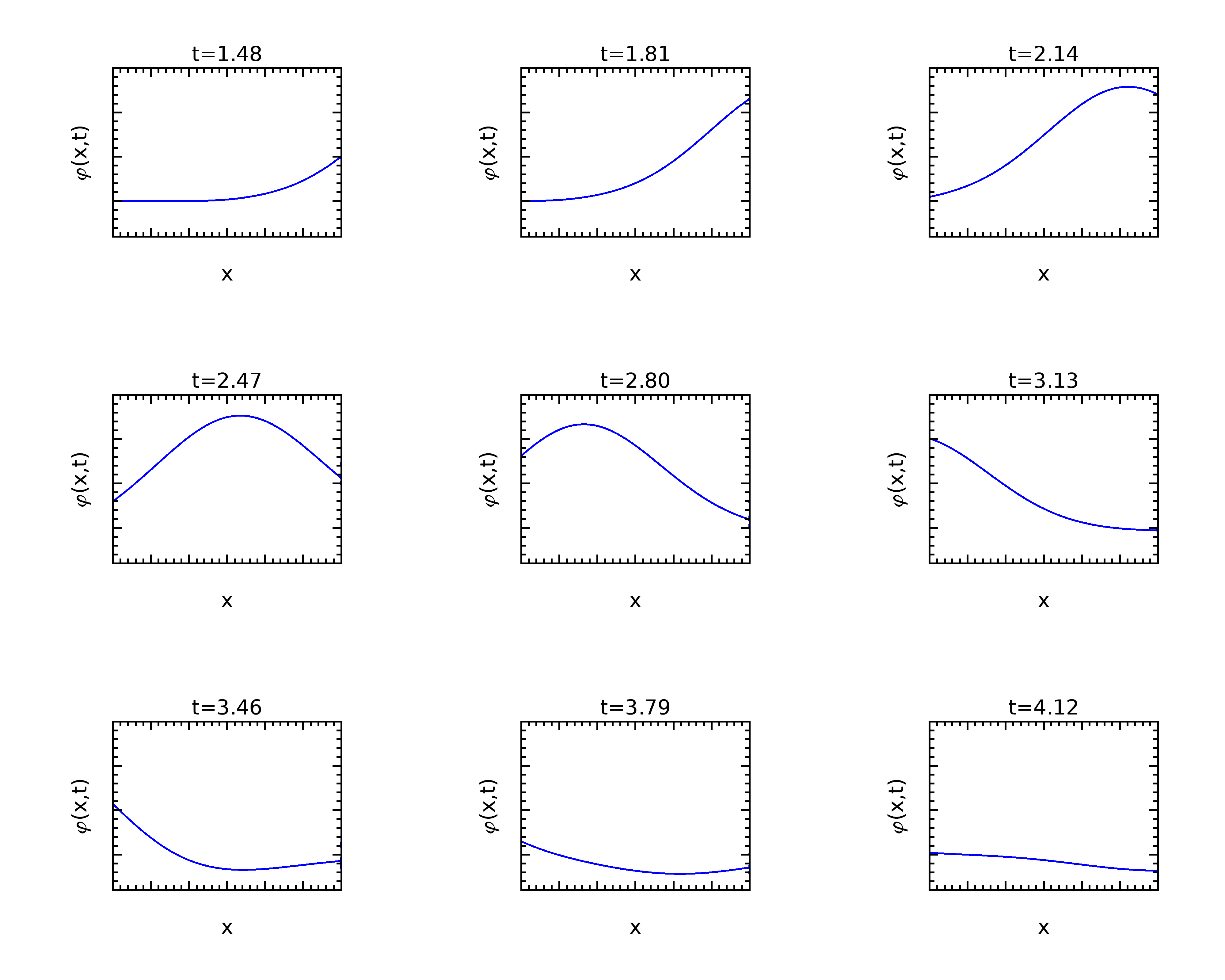}
\caption{A numerical solution of the EOS formulation for model 2 generated by an external source. The parameter values used are  $a_0=0.0, a_1=3.0, N=1600,$ $\alpha=-1.0,$ $\beta=0.3$, $\gamma=8.0$, $\mu=2.0,$ $\nu=2.0,$ $\mu_0=1.0$, $\nu_0=1.0, A=1.0, \alpha_1=36, \beta_1=4, 
t_s=1.0, x_o=4.0$.}
\label{ExampleRunModel2}
\end{figure}

Figure \ref{ExampleRunModel2} shows scattering of a wave generated by an external source calculated from our numerical  implementation of the EOS formulation for model 2. The source we used is given by
\begin{equation*}
j_s=Ae^{-\alpha_1(x-x_o)^2-\beta_1(t-ts)^2}.
\end{equation*}

\section*{Acknowledgements}
All three authors are thankful for  support from the department of mathematics and statistics at UIT the Arctic University of Norway. The third author is also thankful for academic and material support  from  the Arizona Center for Mathematical Sciences at the University of Arizona and for the support from the Air Force Office for Scientific Research under grant \# FA9550-16-1-0088.

\begin{appendices}
\section{Stability of the numerical schemes for model 1 and model 2}
As mentioned in the main text, we don't expect the two 1D models in this paper to be representative for stability issues pertaining to numerical implementation to EOS formulations in general. However there is an issue that is worth discussing here. From the EOS formulation of model 1 one might expect that there would be severe stability issues associated with any numerical approximation. The reason is that the basic equation for the field inside the domain $(a_0,a_1)$
\begin{align}
\varphi_{t}  &  =c_1 \varphi_{x},\label{Model1Appendix}
\end{align}
uncoupled for simplicity from the internally generated current density $j$, can only satisfy the boundary condition at the right boundary $a_1$  induced by the external source. This is because equation (\ref{Model1Appendix})  is of order one in space derivatives and consequently one can not impose any additional boundary condition at $a_0$ that is independent of the one imposed at $a_1$. The EOS formulation evades this problem by imposing, in this simplified setting, the boundary condition 
\begin{align}
\varphi_1(a_{0},t)=\theta(a_{0}-a_{1}+c_1(t-t_{0}))\varphi(a_{1},t-\frac{a_{1}-a_{0}}
{c_1}),\label{A17}
\end{align}
which depends on the boundary condition at $a_1$ in exactly the way it needs in order for a solution, for the boundary value problem for (\ref{Model1Appendix}) on the interval $(a_0,a_1)$,  to exist. However, this existence seems precarious, if we miss the right value by even a small amount in a numerical scheme, are we not then solving a boundary value problem for (\ref{Model1Appendix}) where the two boundary conditions are {\it not} related in the right way, and is there not a danger that this non-existence will manifest itself in a numerical instability? In fact, could it be that the restricted domain of stability of the EOS formulation,  as noted in the main text, is a result of the very particular, in general delay-boundary conditions imposed as a consequence of the EOS formulation? If this was true it would be important because such delayed boundary conditions are a general feature of EOS formulations. We will however now show that the restricted domain of stability for the 1D models discussed in this paper are in fact caused by the nonuniformity rather than the delayed type boundary conditions. 

For this purpose we introduce a family of grids of the interval $(a_0,a_1)$ that are parametrized by $\epsilon$. The grid is uniform for $\epsilon=0$ and is equal to the nonuniform grid we used for our numerical implementations for model 1 and 2 when $\epsilon=1$.
\begin{equation*} 
x_i=a_0+(i+1-0.5 \epsilon)\Delta x ,\  i=0, 1, \cdots, N-1, \,\\
\end{equation*}
where  $\epsilon \in [0,1]$ and 
\begin{equation*}
\Delta x  =\frac{N+\epsilon}{N (N+1)}(a_1-a_0). \, \\
\end{equation*}
 In order to derive a finite difference scheme for (\ref{Model1Appendix}), using the Lax-Wendroff approach like in the main text, we need to impose some boundary conditions. In the end these conditions do not influence the stability of the scheme and we therefore for simplicity  impose fixed boundary conditions. Given this the numerical scheme takes the form
\begin{equation}
{\bf U}_{n+1}=M_1{\bf U}_n+{\bf b},\label{m1}
\end{equation}
where ${\bf  U}=(\bf  \varphi)$
is a N vector, $M_1$ is a matrix of order $N \times N$ given by 
\newcommand\scalemath[2]{\scalebox{#1}{\mbox{\ensuremath{\displaystyle #2}}}}

\begin{align*}
M_1=\left[ 
 \scalemath{0.8}{
\begin{matrix}
\eta_1+c_1\, \eta_2  &\gamma_1+c_1\, \gamma_2 & 0  &0 &0&\ldots&0\\
\kappa_1-c_1\, \kappa_2  & \chi  &\kappa_1+c_1\, \kappa_2&0 &0&\ldots&0\\
0 &\kappa_1-c_1\, \kappa_2  &  \chi& \kappa_1+c_1\, \kappa_2 &0&\ldots &0\\
\vdots&\vdots&\vdots&\vdots&\vdots&\vdots&\vdots\\
0 &\ldots&0&0&\kappa_1-c_1\, \kappa_2  & \chi &\kappa_1+c_1\, \kappa_2\\
0 &\ldots&0&0&0& \gamma_3-c_1\, \gamma_4   &\eta_3-c_1\, \eta_4 
\end{matrix}
}\right] \, \\
\end{align*}
where the entries of the matrix depend on the discrete grid but not on the boundary conditions 
and where ${\bf b}$ is determined by the boundary values.\\

  Let us look for a constant solution to (\ref{m1}), $ {\bf U}={\bf U}^*$. For ${\bf U}^*$ to be a solution, we must have 

\begin{align}
{\bf U}^*=M_1&{\bf U}^*+{\bf b},\nonumber\\
&  \Updownarrow\nonumber\\
(M_1-I)&{\bf U}^*    ={\bf b}\label{m2},
\end{align}
where  $I$ is identity matrix of order $N\times N$.
In order to have a unique solution for (\ref{m2}), $\lambda=1$ must not be an eigenvalue for $M$. Thus, the unique solution will be given by
\begin{equation*}
{\bf U}^*=(M_1-I)^{-1}{\bf b}.
\end{equation*}
Define now ${\bf y}_n$ by 
\begin{equation*}
{\bf U}_n={\bf y}_n+{\bf U}^*.
\end{equation*}

\begin{align}
{\bf y}_{n+1}+{\bf U}^*&=M_1({\bf y}_{n}+{\bf U}^*)+{\bf b},\nonumber\\
&  \Updownarrow\nonumber\\
{\bf y}_{n+1}&=M_1{\bf y}_{n}\label{m3}.
\end{align}
The matrix $M_1$ is not symmetric, but numerical investigations show that it  in general has $N$ different eigenvalues,  $\lambda_i,  i=1,2,\cdots, N $ . Then the corresponding eigenvectors, ${\bf y}_i,$ are then  independent and form a basis for ${\rm I\!R}$. Let now ${\bf y}_0 \in {\rm I\!R}$ be an initial value for (\ref{m3}). Then we have 

\begin{align*}
{\bf y}_0&=\sum_i d_i{\bf y}_i,\nonumber\\
&  \Updownarrow\nonumber\\
{\bf y}_n&=M_1^n {\bf y}_0=\sum_i d_i\lambda_i^n{\bf y}_i.
\end{align*}
We can see that if there exists any eigenvalue that is located outside the unit circle, then $\|{{\bf y}_n}\| \rightarrow \infty$. So in order to get a stable numerical solution to model (\ref{Model1Appendix}), the eigenvalues of $M_1$ must satisfy
\begin{equation*}
\max_{i}\ | \lambda_i \ | <1.
\end{equation*}
Using this result we find that the stability domain as a function of $\epsilon$  is of the form
\begin{align}
\tau_1(\epsilon)\frac{\Delta x}{c_1}<\Delta t<\tau_2(\epsilon)\frac{\Delta x}{c_1}\label{stabilitydomain}.
\end{align}

\begin{figure}[h!]
\centering
\includegraphics[height=10cm, width=10cm]{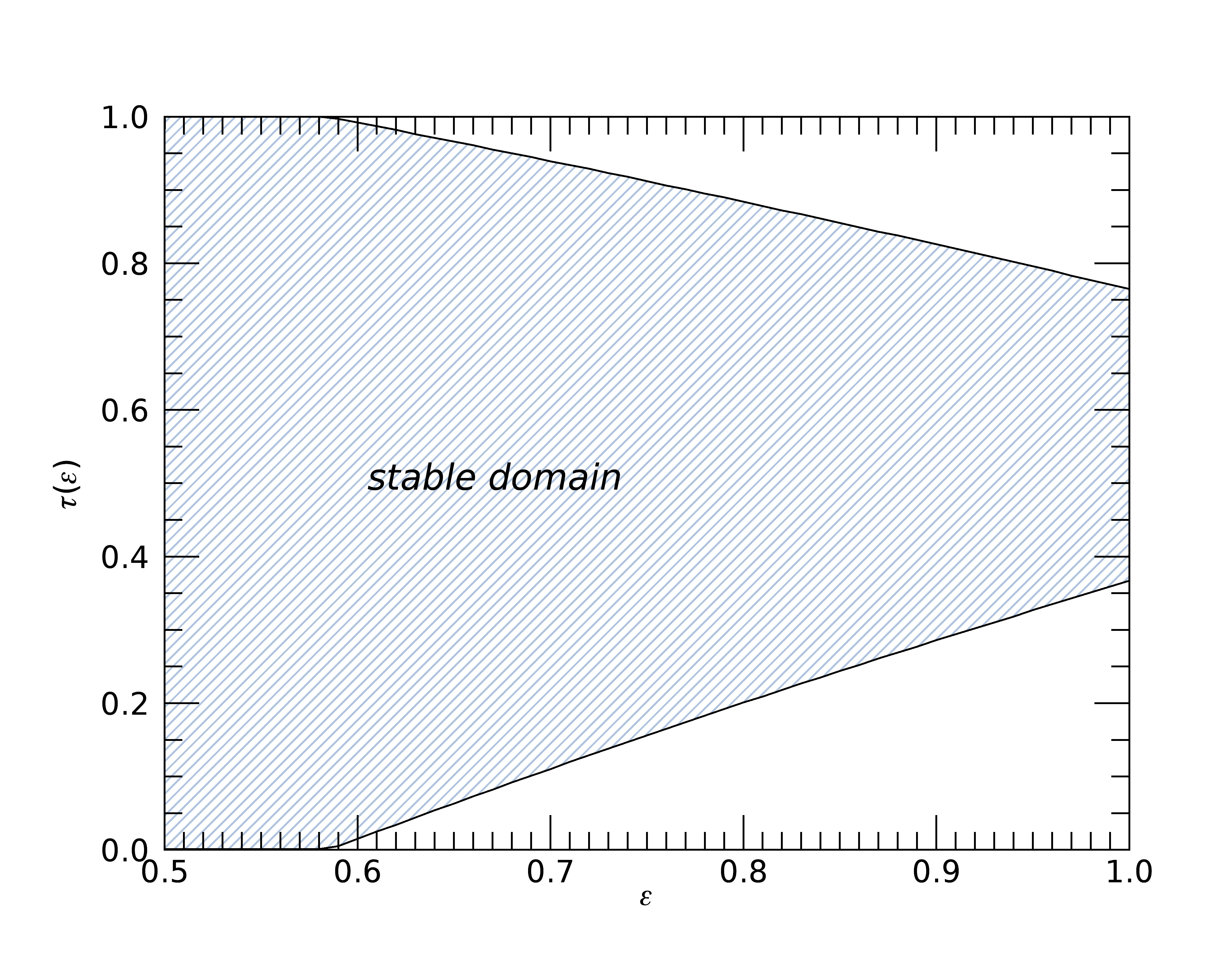}
\caption{ The stability domain for the EOS formulation of (\ref{Model1Appendix})}
\label{Stabledomain}
\end{figure}

It is evident from figure (figure \ref{Stabledomain})  that the restriction on the stability domain for the EOS formulation of (\ref{Model1Appendix}),  as compared to the Lax-Wendroff scheme for the case of free space propagation, is caused by the introduction of a nonuniform grid for the EOS formulation.

For model 2 we find exactly the same stability domain as illustrated in figure (figure \ref{Stabledomain}) for model 1. That there should be some relation between the stability of these two models is perhaps not very surprising at the level of PDEs. After all, if we decouple the fields in model 2 from the current, the resulting system is equivalent to the wave equation and solutions of that equation are sums of left and right going waves of the type described by Model 1. However, at the level of numerical schemes the coinciding  of the stability domains for the two models is somewhat less obvious. Note that we can write the matrix $M_1$, determining the stability for model 1, in the form
$$
M_1=m_1+c\,m_1,
$$
where $m_1$ and $m_2$ are $N\times N$ matrices given by,
$$
m_1=
\left[ \begin{matrix}
\eta_1  &\gamma_1 & 0  &0 &0&\ldots&0\\
\kappa_1  & \chi  &\kappa_1&0 &0&\ldots&0\\
0 &\kappa_1 &  \chi& \kappa_1 &0&\ldots &0\\
\vdots&\vdots&\vdots&\vdots&\vdots&\vdots&\vdots\\
0 &\ldots&0&0&\kappa_1  & \chi &\kappa_1\\
0 &\ldots&0&0&0& \gamma_3  &\eta_3 
\end{matrix}\right],
$$
$$
m_2=
\left[ \begin{matrix}
\eta_2  &\, \gamma_2 & 0  &0 &0&\ldots&0\\
-\, \kappa_2  & 0 &\, \kappa_2&0 &0&\ldots&0\\
0 &-\, \kappa_2  &  0& \kappa_2 &0&\ldots &0\\
\vdots&\vdots&\vdots&\vdots&\vdots&\vdots&\vdots\\
0 &\ldots&0&0&-\, \kappa_2  &0 &\, \kappa_2\\
0 &\ldots&0&0&0& -\, \gamma_4   &-\, \eta_4 
\end{matrix}\right].
$$
Given this, the $2 N\times 2 N$ matrix determining the stability of model 2 is given by
$$
M_2=
\left[ \begin{matrix}
m_1&\mu_1 \, m_2\\
-\nu_1 \, m_2 &m_1
\end{matrix}\right].
$$
The matrix $M_2$ clearly has a block structure and the same blocks give a linear decomposition of $M_1$ into a sum of two terms. However, we were not able use these commonalities between $M_1$ and $M_2$ to explain the fact that model 1 and model 2 have, not approximately, but exactly the same domain of stability as far as we can determine.
Note that the occurence of a stability domain like (\ref{stabilitydomain}) might be a universal feature of EOS formulations. We have for example found a stability domain of this type in the EOS formulation of 3D Maxwell's equations. There, however, it is clear that the delay boundary conditions is at least in part responsible for the width of the stability domain.

\section{Numerical implementation of the EOS formulation for model 2}

The numerical implementation of model 2 contains the same elements as the ones we introduced  for model 1. Thus we first  define a nonuniform space grid  inside  the scattering object, $(a_0,a_1)$,
\begin{equation}
x_i=a_0+(i+0.5)\Delta x,\  i=0, 1, \cdots, N-1\label{TwoFieldGrid},
\end{equation}
where $\Delta x =\frac{a_1-a_0}{N}$. The grid points (\ref{TwoFieldGrid}) are the internal nodes for model 2. We also introduce the discrete time grid 
\begin{equation*}
t^n=n\Delta t,\  n=0, 1, \cdots.  
\end{equation*}
The values of the  parameter $\Delta t$ will of course, like for model 1,  be bounded by the requirement of stability for the scheme.

We apply the Lax-Wendroff method to the first three equations of (\ref{B1.1}) and the modified Euler's method to the last equation of (\ref{B1.1}). For interval $(a_0, a_1)$, the numerical iteration  can be written as
\begin{align}
\varphi_i^{n+1}=&\varphi_i^n+\Delta t \ (\mu_1\frac{\partial \psi}{\partial x}+j )_i^n+\frac{1}{2} (\Delta t)^2( \mu_1 \nu_1 \frac{\partial^2 \varphi}{\partial x^2}+f)_i^n,\nonumber\\
\psi_i^{n+1}=&\psi_i^n+\Delta t \ (\nu_1 \frac{\partial \varphi}{\partial x})_i^n+\frac{1}{2} (\Delta t)^2( \mu_1 \nu_1 \frac{\partial^2 \psi}{\partial x^2}+\nu_1 \frac{\partial j}{\partial x})_i^n,\nonumber\\
\rho_i^{n+1}=&\rho_i^n+\Delta t \ ( - \frac{\partial j}{\partial x} )_i^n+\frac{1}{2} (\Delta t)^2( - \frac{\partial f}{\partial x})_i^n,\nonumber\\
\bar j_i^{n+1}=&j_i^n+\Delta t \ f_i^n,\nonumber\\
j_i^{n+1}=&\frac{1}{2}(j_i^n+\bar j_i^{n}+\Delta t\  f(\rho_i^{n+1},\varphi_i^{n+1},\bar j_i^{n+1}))\nonumber,
\end{align}
where $f=(\alpha -\beta \rho)\varphi-\gamma j$. The finite difference approximations for the fields and the current density at all internal nodes,except the two nodes closest to the boundary points $a_0$ and $a_1$,  are given by the standard expressions
\begin{align*}
(\frac{\partial \phi}{\partial x})_i^n&=\frac{\phi_{i+1}^n-\phi_{i-1}^n}{2 \Delta x},\nonumber\\
(\frac{\partial^2 \phi}{\partial x^2})_i^n&=\frac{\phi_{i+1}^n-2 \phi_i^n+\phi_{i-1}^n}{(\Delta x)^2},\ \phi =\varphi, \psi,  j ,
\end{align*}
for $ i=1,2,\cdots ,N-2,$. For the two internal nodes closest to the boundary points, we need to use alternative difference rules because the grid is nonuniform in the domain around these nodes
\begin{align*}
&(\frac{\partial \phi}{\partial x})_0^n=-\frac{1}{3\cdot \Delta x}(4\phi_{a_0}^n-3\phi_0^n-\phi_1^n),\nonumber\\
&(\frac{\partial^2 \phi}{\partial x^2})_0^n=\frac{4}{3\cdot (\Delta x)^2}(2\phi_{a_0}^n-3\phi_0^n+\phi_1^n),\nonumber\\
&(\frac{\partial \phi}{\partial x})_{N-1}^n=\frac{1}{3\cdot \Delta x}(4\phi_{a_1}^n-3\phi_{N-1}^n-\phi_{N-2}^n),\nonumber\\
&(\frac{\partial^2 \phi}{\partial x^2})_{N-1}^n=\frac{4}{3\cdot (\Delta x)^2}(2\phi_{a_1}^n-3\phi_{N-1}^n+\phi_{N-2}^n),\nonumber\\
&(\frac{\partial j}{\partial x})_0^n=\frac{1}{2\Delta x}(4j_1^n-3j_0^n-j_2^n),\nonumber\\
&(\frac{\partial j}{\partial x})_{N-1}^n=-\frac{1}{2 \Delta x}(4j_{N-2}^n-3j_{N-1}^n-j_{N-3}^n),
\end{align*}
where $\phi=\varphi, \psi$.
The discretization of the boundary update rules (\ref{B38.1}) and  (\ref{B38.2}) are
\begin{align*}
&\begin{pmatrix}
c_1+c_0 &\mu_1-\mu_0\\
\nu_1-\nu_0&c_1+c_0
\end{pmatrix}
\dbinom{\varphi}{\psi}(a_0, t_{n+1})\\
\nonumber\\
&=\frac{\Delta x}{c_1}
\begin{pmatrix}
c_1 &\mu_1\\
\nu_1&c_1
\end{pmatrix}\sum_{i=0}^{N-1} \theta(t_{n+1}-t_0-\frac{x_i-a_0}{c_1})\dbinom{j}{0}(x_i,t_{n+1}-\frac{x_i-a_0}{c_1})\\
\nonumber\\
&+
\begin{pmatrix}
c_1&\mu_1\\
\nu_1&c_1
\end{pmatrix}\theta(t_{n+1}-t_0-\frac{a_1-a_0}{c_1})
\dbinom{\varphi}{\psi}_-(a_1,t_{n+1}-\frac{a_1-a_0}{c_1})
\end{align*}
\begin{align*}
&\begin{pmatrix}
c_1+c_0 &\mu_0-\mu_1\\
\nu_0-\nu_1&c_1+c_0
\end{pmatrix}
\dbinom{\varphi}{\psi}(a_1, t_{n+1})=\\
\nonumber\\
&\frac{\Delta x}{c_1}
\begin{pmatrix}
c_1 &-\mu_1\\
\nonumber\\
-\nu_1&c_1
\end{pmatrix}\sum_{i=0}^{N-1} \theta(t_{n+1}-t_0-\frac{a_1-x_i}{c_1})\dbinom{j}{0}(x_i,t_{n+1}-\frac{a_1-x_i}{c_1})\\
\nonumber\\
&-
\begin{pmatrix}
-c_1&\mu_1\\
\nu_1&-c_1
\end{pmatrix}\theta(t_{n+1}-t_0-\frac{a_1-a_0}{c_1})
\dbinom{\varphi}{\psi}_-(a_0,t_{n+1}-\frac{a_1-a_0}{c_1}) \\
\nonumber\\
&+2c_{0}\left(
\begin{array}
[c]{c}%
\varphi_{i}\\
\psi_{i}%
\end{array}
\right)  (a_{1},t_{n+1}).
\end{align*}
where $\left(
\begin{array}
[c]{c}%
\varphi_{i}\\
\psi_{i}%
\end{array}
\right)  (a_{1},t_{n+1})$ are determined by the external source. 

\section{Artificial source test of the EOS formulation for model2}

The source extended model 2,  is given by
\begin{align*}
\varphi_t&=\mu_1 \psi_x +j+g_1\nonumber,\\
\psi_t &=\nu_1 \varphi_x +g_2\nonumber,\\
\rho_t & =-j_x+g_3\nonumber,\\
j_t &=(\alpha-\beta \rho)\varphi-\gamma j+g_4.
\end{align*}
For the source extended model 2, any given set of functions $\varphi_0,$ $\psi_0,$ $j_0$ and $\rho_0$ is a solution if the sources are chosen to be 
\begin{align*}
g_{01}&=(\varphi_0)_t-\mu_1 (\psi_0)_x -j_0\nonumber,\\
g_{02}&=(\psi_0)_t -\nu_1 (\varphi_0)_x \nonumber,\\
g_{03}&=(\rho_0)_t  +(j_0)_x\nonumber,\\
g_{04}&=(j_0)_t -(\alpha-\beta \rho_0)\varphi_0+\gamma j_0.
\end{align*}
The boundary update rule for the source extended model 2 is changed into

\begin{align*}
&\left(
\begin{array}
[c]{cc}%
c_1+c_{0} & \mu_1-\mu_0\\
\nu_1-\nu_0 & c_1+c_{0}%
\end{array}
\right)  \left(
\begin{array}
[c]{c}%
\varphi_{1}\\
\psi_{1}%
\end{array}
\right)  (a_{0},t)=\label{B37}\nonumber\\
\nonumber\\
&\frac{1}{c_1}\left(
\begin{array}
[c]{cc}%
c_1 & \mu_1\\
\nu_1 & c_1%
\end{array}
\right)
{\displaystyle\int_{a_{0}}^{a_{1}}}
dx^{\prime}\theta(c_1(t-t_{0})-(x^{\prime}-a_{0}))\binom{j+g_1}{g_2}(x^{\prime
},t-\frac{x^{\prime}-a_{0}}{c_1})\nonumber\\
\nonumber\\
&+\theta(c_1(t-t_{0})-(a_{1}-a_{0}))\left(
\begin{array}
[c]{cc}%
c_1 & \mu_1\\
\nu_1 & c_1%
\end{array}
\right)  \left(
\begin{array}
[c]{c}%
\varphi_{1}\\
\psi_{1}%
\end{array}
\right)  (a_{1},t-\frac{a_{1}-a_{0}}{c_1}),\\
\end{align*}
\begin{align*}
&\left(
\begin{array}
[c]{cc}%
c_{0}+c_1 & \mu_0-\mu_1\\
\nu_0-\nu_1 & c_1+c_{0}%
\end{array}
\right)  \left(
\begin{array}
[c]{c}%
\varphi_{1}\\
\psi_{1}%
\end{array}
\right)  (a_{1},t)=\nonumber\\
\nonumber\\
&\frac{1}{c_1}\left(
\begin{array}
[c]{cc}%
c_1 & -\mu_1\\
-\nu_1 & c_1%
\end{array}
\right)
{\displaystyle\int_{a_{0}}^{a_{1}}}
dx^{\prime}\theta(c_1(t-t_{0})-(a_{1}-x^{\prime}))\binom{j+g_1}{g_2}(x^{\prime
},t-\frac{a_{1}-x^{\prime}}{c_1})\nonumber\\
\nonumber\\
&-\theta(c_1(t-t_{0})-(a_{1}-a_{0}))\left(
\begin{array}
[c]{cc}%
-c_1 & \mu_1\\
\nu_1 & -c_1%
\end{array}
\right)  \left(
\begin{array}
[c]{c}%
\varphi_{1}\\
\psi_{1}%
\end{array}
\right)  (a_{0},t-\frac{a_{1}-a_{0}}{c_1})\nonumber\\
\nonumber\\
&+2c_{0}\left(
\begin{array}
[c]{c}%
\varphi_{i}\\
\psi_{i}%
\end{array}
\right)  (a_{1},t). 
\end{align*}

\end{appendices}

\bibliographystyle{plain}
\bibliography{paper1Aihua}

\end{document}